\ifdefined\journalstyle
   \documentclass{mcom-l}
\else
  \documentclass[11pt]{article} 
  \def\arxivstyle{1}
  \usepackage{theorem}
  \newenvironment{proof}{\begin{trivlist}
  \item[\hskip\labelsep{\it Proof.}]}{$\hfill\Box$\end{trivlist}}
\fi

\usepackage{pdfsync}
\usepackage[hidelinks]{hyperref}

\usepackage{latexsym,amsfonts,amsmath,amssymb,url}
\usepackage{graphicx}
\usepackage{subcaption} 

\usepackage[usenames,dvipsnames,svgnames]{xcolor}

\newtheorem{theorem}{Theorem}[section]       
\newtheorem{lemma}[theorem]{Lemma}

\numberwithin{equation}{section}             

\newenvironment{proofof}[1]{\begin{trivlist}
		\item[\hskip\labelsep{\bf Proof of {#1}.}]}{$\hfill\Box$\end{trivlist}}

\newcommand{\indx}{\mathcal{I}}
\newcommand{\satop}[2]{\stackrel{\scriptstyle{#1}}{\scriptstyle{#2}}}

\newcommand{\bsalpha}{{\boldsymbol{\alpha}}}
\newcommand{\bsbeta}{{\boldsymbol{\beta}}}

\newcommand{\bsb}{{\boldsymbol{b}}}

\newcommand{\bse}{{\boldsymbol{e}}}

\newcommand{\bsgamma}{{\boldsymbol{\gamma}}}
\newcommand{\bslam}{{\boldsymbol{\lambda}}}

\newcommand{\bsnu}{{\boldsymbol{\nu}}}

\newcommand{\bst}{{\boldsymbol{t}}}

\newcommand{\bsu}{{\boldsymbol{u}}}
\newcommand{\bsv}{{\boldsymbol{v}}}
\newcommand{\bsx}{{\boldsymbol{x}}}

\newcommand{\bsy}{{\boldsymbol{y}}}
\newcommand{\bsz}{{\boldsymbol{z}}}
\newcommand{\bsw}{{\boldsymbol{w}}}

\newcommand{\bsm}{{\boldsymbol{m}}}

\newcommand{\bszero}{{\boldsymbol{0}}}

\newcommand{\rd}{{\mathrm{d}}}

\newcommand{\bbA}{{\mathbb{A}}}
\newcommand{\bbB}{{\mathbb{B}}}
\newcommand{\bbN}{{\mathbb{N}}}

\newcommand{\bbR}{{\mathbb{R}}}

\newcommand{\calE}{{\mathcal{E}}}
\newcommand{\calI}{{\mathcal{I}}}
\newcommand{\calJ}{{\mathcal{J}}}
\newcommand{\calK}{{\mathcal{K}}}

\newcommand{\calG}{{\mathcal{G}}}

\newcommand{\calO}{{\mathcal{O}}}

\newcommand{\calR}{{\mathcal{R}}}
\newcommand{\calS}{{\mathcal{S}}}
\newcommand{\calW}{{\mathcal{W}}}

\newcommand{\mask}[1]{}

\newcommand{\setu}{{\mathfrak{u}}}

\DeclareMathOperator*{\argmin}{\mathrm{arg\,min}}

\definecolor{darkred}{RGB}{139,0,0}
\definecolor{darkgreen}{RGB}{0,100,0}
\definecolor{darkmagenta}{RGB}{180,0,180}
\definecolor{darkblue}{RGB}{0,0,190}
\definecolor{darkorange}{RGB}{180,60,0}

\begin{document}

\date{March 2026}

\ifdefined\journalstyle
  \title[Regularity and tailored regularization of Deep Neural Networks]
  {Regularity and tailored regularization of Deep Neural Networks, with application to
   parametric PDEs in uncertainty quantification}
  \author{Alexander Keller}
  \address{NVIDIA, Berlin, Germany}
  \email{akeller@nvidia.com}
  \author{Frances Y. Kuo}
  \address{School of Mathematics and Statistics, UNSW Sydney, Sydney, Australia}
  \email{f.kuo@unsw.edu.au} 
  \author{Dirk Nuyens} 
  \address{Department of Computer Science, KU Leuven, Leuven, Belgium}
  \email{dirk.nuyens@kuleuven.be}
  \author{Ian H. Sloan}  
  \address{School of Mathematics and Statistics, UNSW Sydney, Sydney, Australia}
  \email{i.sloan@unsw.edu.au}
\else
  \title{Regularity and tailored regularization of Deep Neural Networks, with application to
   parametric PDEs in uncertainty quantification}
  \author{Alexander Keller, Frances Y. Kuo, Dirk Nuyens, and Ian H. Sloan}
\fi

\ifdefined\arxivstyle
  \maketitle  
\fi

\begin{abstract}
In this paper we consider Deep Neural Networks (DNNs) with a smooth
activation function as surrogates for high-dimensional functions that are
somewhat smooth but costly to evaluate. We consider the standard
(non-periodic) DNNs as well as propose a new model of periodic DNNs which
are especially suited for a class of periodic target functions when
Quasi-Monte Carlo lattice points are used as training points. 
The primary contribution of this paper is the derivation of explicit bounds
for all mixed derivatives of DNNs with respect to their input parameters.
The bounds depend on the neural network parameters as well as the choice of activation function, with explicit constants.
These bounds are fully general and remain independent of both the
target function and the training data. 
By imposing restrictions on the network parameters to
match the regularity features of the target functions, we prove that DNNs
with $N$ tailor-constructed lattice training points can achieve the
generalization error (or $L_2$ approximation error) bound ${\tt tol} +
\calO(N^{-r/2})$, where ${\tt tol}\in (0,1)$ is the tolerance achieved by
the training error in practice, and $r = 1/p^*$, with $p^*$ being the
\emph{summability exponent} of a sequence that characterises the decay of
the input variables in the target functions, and with the implied constant
independent of the dimensionality of the input data. We apply our analysis
to popular models of parametric elliptic PDEs in uncertainty
quantification. In our numerical experiments, we restrict the network
parameters during training by adding tailored regularization terms, and we
show that for an algebraic equation mimicking the parametric PDE problems
the DNNs trained with tailored regularization perform significantly
better.
\end{abstract}

\ifdefined\journalstyle
  \maketitle 
\else

  \smallskip\noindent
  \textbf{MSC2020 Subject Classification:} 68T07, 65D32, 65D40, 35R60
\fi


\section{Introduction}

We seek an efficient way to approximate Data-to-Observables maps (DtOs)
for a large class of physics-based problems of the general form
\begin{equation} \label{eq:DtOs}
\mbox{Given data $\bsy\in Y$, compute observable $G(\bsy)\in Z$},
\end{equation}
where $\bsy\in Y \subset\bbR^s$ is a vector of $s$ input parameters to a
parametric PDE, and $G(\bsy)\in Z = \bbR^{N_{\mathrm{obs}}}$ is a vector
of observables obtained from the solution to a PDE which is a somewhat
smooth function of $\bsy$. The observable may be the average of the PDE
solution over a subset of the domain, in which case $N_{\mathrm{obs}}=1$;
or it may be the vector of values of the solution $u$ at the points of a
finite element mesh, in which case $N_{\mathrm{obs}}$ may be the number of
interior finite element mesh points.

Problems of this kind can be extremely challenging, especially when a high
cost for one PDE solution is combined with hundreds of uncertain input
parameters. Recent reviews of parametric PDEs \cite{ABW22,CD15} both refer
to the ``curse of dimensionality'' for large-scale computations of this
kind.
Thus a major theme of recent international research for parametric PDEs
has been the search for fast surrogates. This research comes under many
different headings, including best $N$-term approximation \cite{CDS10} and
reduced order modelling \cite{GWW23,GW21}. Multi-level
methods~\cite{Gil15} and multi-fidelity methods \cite{PWG18} play
a big role, too.

This paper analyzes a popular approach to approximate $G(\bsy)$, namely, via a \emph{Deep Neural
Network} (DNN) \cite{GBC16,LBH15}. We interpret DNN mathematically as a
nonlinear algorithm for multivariate function approximation, and we aim to
derive conditions on the DNN to obtain favorable theoretical error bounds.
Indeed, in the search for faster and more reliable methods for parametric
PDEs, use is increasingly being made of machine learning
\cite{GK22,MT18,SB14}, under the guiding principle that the neural network
should be ``physics informed'' -- the physics underlying the PDE system
should be respected \cite{MM23,RPK19}. We follow the papers
\cite{LMR20,LMRS21,MR21} which used \emph{Quasi-Monte Carlo} (QMC) methods
to ``train'' the DNN, and in assuming that the training points can be
chosen freely by the user, as for example in synthetic data generation
(SDG) by simulation, e.g., \cite{IEQandML,MESK22,MRNK21}. 
However, there are two crucial differences from
\cite{LMR20,LMRS21,MR21}. Firstly, in addition to the standard
(non-periodic) DNNs, we propose a new model of periodic DNNs and we use
QMC points that are adapted to the periodicity. Secondly, we obtain
explicit regularity bounds on both models of DNNs, and we impose
restrictions on the neural network parameters to match the regularity
features of the target functions. 

To the best of our knowledge, this is the first paper to establish parametric regularity bounds on DNNs with explicit constants that depend on the network parameters. These bounds are entirely general, remaining independent of the target function and the selection of training points. While we use these bounds to develop a rigorously justified QMC theory for training, these general bounds have broader applicability and may be leveraged in other settings, such as sparse grids and best-$n$ term approximations.

\subsection{The framework}

We consider a standard (``feed-forward'') DNN of ``depth''~$L$ as an
approximation of $G$ in \eqref{eq:DtOs} of the form
\begin{align} \label{eq:DNN-np}
  G_\theta^{[L]}(\bsy)
  := W_L\, \sigma (\cdots W_2\, \sigma (W_1\, \sigma (W_0\,\bsy + \bsv_0) + \bsv_1) + \bsv_2 \cdots) + \bsv_L,
  \quad
  \bsy\in Y,
\end{align}
where for each $\ell=0,\ldots,L$ there is an affine transformation with a
``weight'' matrix $W_\ell \in \bbR^{d_{\ell+1} \times d_\ell}$ and a
``bias'' vector $\bsv_\ell\in\bbR^{d_{\ell+1}\times 1}$, and there is a
scalar nonlinear ``activation function'' $\sigma$ to be applied
component-wise between the ``layers''. Layer $0$ is the ``input'' layer,
with $d_0 = s$ the dimensionality of the input vector $\bsy$. Layer $L+1$
is the ``output'' layer, with dimensionality $d_{L+1} = N_{\rm obs}$ being
the number of observables. The $L$ layers in between are the ``hidden''
layers. The values of $d_0,\ldots,d_{L+1}$ are referred to as
the ``widths''. This describes a general fully-connected neural network
(multi-layer perceptron), with $
 N_{\rm DNN} := \sum_{\ell=0}^L (d_{\ell+1}\times d_\ell) + \sum_{\ell=0}^L
 d_{\ell+1}
$
network parameters to be ``learned'' --- these are collectively denoted by
\[
  \theta := (W_\ell,\bsv_\ell)_{\ell=0}^L,
\]
and the set of all permissible network parameters is denoted by $\Theta$.
The depth $L$, the widths $d_\ell$, and other parameters used in the
``training'' phase (see below) are known together as the
``hyperparameters''.
Other specific types of neural networks are defined by imposing
constraints on the structure of the matrices. For example, some networks
assume that all matrices on different layers are of the same size/depth,
and some convolutional neural networks assume that the matrices are
circulant so that the fast Fourier transform can be used to speed up the
calculation.
Popular activation functions include e.g., rectified linear unit (ReLU)
$\sigma(x) = \max\{x,0\}$, logistic sigmoid $\sigma(x) = 1/(1+e^{-x})$,
hyperbolic tangent $\sigma(x) = \tanh(x) = (e^x-e^{-x})/(e^x+e^{-x})$, and
swish $\sigma(x) = x/(1+e^{-x})$.

We have chosen to index $W_\ell$ and $\bsv_\ell$ starting unconventionally
from $\ell=0$ for the benefit of our later presentation, with a minor
pecularity that $\bsv_\ell$ has dimension $d_{\ell+1}$. Our definition of
depth $L$ is the number of hidden layers, which is also the number of
times the activation function $\sigma$ is applied. Each hidden layer
$\ell$ is given by $\sigma(W_{\ell-1}\, \bsy + \bsv_{\ell-1})$ for
$\bsy\in\bbR^{d_{\ell-1}\times 1}$. An interesting modification to
\eqref{eq:DNN-np}, which we leave for future investigation, is to add an
identity ``skip link''~\cite{ResNet1,ResNet2,StableNN} so that layer
$\ell$ is given instead by $\bsy + \sigma(W_{\ell-1}\, \bsy +
\bsv_{\ell-1})$ for $\bsy\in\bbR^{d_{\ell-1}\times 1}$ when $d_{\ell}=d_{\ell-1}$.
Skip links are known
to be essential for training deep networks, since they help mitigate vanishing and
exploding gradients~\cite{LSTM97}.

A training set is a finite set of points $\bsy_1, \ldots, \bsy_N\in Y$
which will be used together with matching (approximate) values of the
observables $G(\bsy_1), \ldots, G(\bsy_N)$ to ``train'' a DNN, i.e., to
determine suitable network parameters $\theta$ to ensure that
$G_\theta^{[L]}(\bsy)$ is a good approximation to $G(\bsy)$, especially in
points that were not present in the training set. We refer to the pairs
$\{(\bsy_k, G(\bsy_k))\}_{k=1}^N$ as our training data.

We follow \cite{LMR20,LMRS21,MR21} in assuming that the parametric input
data is appropriately scaled to the box
\[
  Y := [0,1]^s
\]
of Lebesgue measure $|Y|=1$, and use deterministic points
$\bsy_1,\ldots,\bsy_N$ as training data for a DNN instead of uniformly distributed random points.
To say that the training set $\bsy_1, \ldots, \bsy_N$ is a set of
QMC points is to say that these points are the points of an equal-weight
quadrature rule that approximates the integral over the unit cube: if $F$
is a continuous function in the unit cube then $\int_{Y} F(\bsy)\,\rd \bsy
\approx\frac{1}{N} \sum_{k=1}^N F(\bsy_k)$.

A DNN is ``trained'' by minimising an appropriate ``loss function''. While
our theory developed in later sections allows for an arbitrary number of observables
$N_{\rm obs}$, for simplicity in the introduction we only illustrate the
case $N_{\rm obs}=1$, and we consider the loss function
\[
 \calJ := \frac{1}{N}\sum_{k=1}^N
 \big( G(\bsy_k) - G_\theta^{[L]}(\bsy_k)\big)^2.
\]
The quality of the DNN is often expressed in terms of the
\emph{generalization error} (or $L_2$ approximation error)
\[
 \calE_G := \bigg(\int_Y \big(G(\bsy) - G_\theta^{[L]}(\bsy)\big)^2\,\rd\bsy\bigg)^{1/2}
= \calE_T + (\calE_G - \calE_T),
\]
where $\calE_T := \calJ^{1/2}$ is the computable \emph{training error},
and $(\calE_G - \calE_T)$ is the so-called \emph{generalization gap} which
is related to the quadrature error for the integrand $F(\bsy) := (G(\bsy)
- G_\theta^{[L]}(\bsy))^2$ and which as a result may be reduced with a
well chosen quadrature point set.
It is common in ``supervised'' physics based learning to choose
independent and identically distributed (i.i.d.)~random points $\bsy_k$ as
the training points. With QMC points instead of random points~\cite{IEQandML}, the
convergence rate for the generalization gap (and hence the generalization
error) may
be improved: the papers \cite{LMR20,MR21} used
low discrepancy sequences such as the Sobol$'$ sequence, while the paper
\cite{LMRS21} used a family of high order QMC rules.

In this paper we are particularly interested in a family of QMC rules
called (rank-$1$) \emph{lattice rules} \cite{SJ94} whose points are given
by the simple formula
\begin{align} \label{eq:lat1}
 \bsy_k := \frac{k\bsz \bmod N}{N},\quad k = 1, \ldots, N,
\end{align}
where $\bsz$ is an integer vector of length $s$ chosen from $\{1, 2,
\ldots,N-1\}^{s}$. For convenience, we sometimes use the index $0$ so
that $\bsy_0 \equiv \bsy_N = \bszero$.

\subsection{Our contribution}

In recent works \cite{HHKKS24,KKKNS22,KKS20,KKS24} in the context of
parametric PDEs in uncertainty quantification, the authors introduced a
periodic model of the random field to exploit the full power of lattice
rules. Inspired by these works, in this paper we propose a periodic
variant of the DNN~\eqref{eq:DNN-np} by introducing an extra sine mapping
applied componentwise to the input parameters: we define
\begin{align} \label{eq:DNN-per}
  G_\theta^{[L]}(\bsy)
  := W_L\, \sigma (\cdots W_2\, \sigma (W_1\, \sigma (W_0\,\sin(2\pi\bsy) + \bsv_0) + \bsv_1) + \bsv_2 \cdots) + \bsv_L,
  \quad \bsy\in Y.
\end{align}
This nonlinear function is then ``one-periodic'' with respect to the input
parameters~$\bsy$, i.e.,
$G_\theta^{[L]}(y_1,\ldots,y_{j-1},0,y_{j+1},\ldots,y_s) =
G_\theta^{[L]}(y_1,\ldots,y_{j-1},1,y_{j+1},\ldots,y_s)$ for all
coordinate indices $j=1,\ldots,s$. The particular form \eqref{eq:DNN-per}
is designed to handle a target function $G(\bsy)$ that is itself a
function of $\sin(2\pi\bsy)$. While it may not be suitable for all periodic integrands in general, the resulting DNN is directly suitable for the quantities of interest from
\cite{HHKKS24,KKKNS22,KKS20,KKS24}. It is not intended to be a universal
choice for all periodic integrands.
The extra sine mapping can be viewed as an extra periodic activation function on the input layer,
sometimes known as ``sinusoidal representation networks'' or ``Sirens'',
see e.g., \cite{AMBLW20}. Note that the activation function $\sigma$ in all subsequent layers is generally non-periodic.
Henceforth in this paper, we shall refer to \eqref{eq:DNN-per} as ``the
periodic DNN'', and to \eqref{eq:DNN-np} as ``the non-periodic DNN''.

Section~\ref{sec:reg} is devoted to regularity bounds. The main result of
this paper is Theorem~\ref{thm:der}, which proves explicit bounds on the
mixed derivatives of both the non-periodic and periodic DNNs. To ensure
that these regularity bounds not only are finite but also do not grow with
increasing dimensionality $s$ or increasing network depth~$L$ or widths,
we need to keep three important factors under control:
\begin{itemize}
\item \textbf{The columns of the matrix $W_0$ need to decay} in vector
    $\infty$-norm. These are characterized by the sequence
    $(\beta_j)_{j\ge 1}$, see \eqref{eq:beta} ahead. This sequence
    moderates the diminishing importance of the successive variables
    in $\bsy$. As the nominal dimension $s$ increases, we need
    $\beta_j \to 0$ as $j\to\infty$. The faster the decay, the
    smaller is the ``effective dimension'' \cite{CMO97}.

\item \textbf{The matrices $W_1,\ldots,W_{L-1}$ need to be bounded} in
    matrix $\infty$-norm. These are characterized by the sequence
    $(R_\ell)_{\ell\ge 1}$, see \eqref{eq:R} ahead. As the network
    increases in depth~$L$ or as we consider even higher order
    derivatives, we prefer $R_\ell\le \rho$ for some small
    constant $\rho>0$ to avoid their product from becoming too large.
    Since the matrix $\infty$-norm is the maximum absolute row sum, if
    the network widths increase then necessarily the matrix elements
    will need to decrease in magnitude. The last matrix $W_L$ does not
    need to be restricted as much as the others. In fact, the last
    matrix needs to account for the scale of the learned function in a
    setting with bounded activation functions. 

\item \textbf{The derivatives of the activation function need to be
    bounded}. These are characterized by the sequence $(A_n)_{n\ge
    1}$, see \eqref{eq:sigma} ahead. For example, the sigmoid
    activation function satisfies $|\sigma^{(n)}(x)|\le n!$ which is
    manageable.
\end{itemize}
Both of the explicit regularity bounds in
Theorem~\ref{thm:der} are new; both depend explicitly on the choice of
activation function. In Theorem~\ref{thm:reg} we state the regularity
bounds in the form that covers three commonly used activation functions, namely, sigmoid, tanh and swish. We remark that these bounds take a very familiar form, similar to those of
regularity bounds of solutions to parametric PDEs in many recent papers.

Section~\ref{sec:err} is devoted to error analysis. Recall that the
generalization gap $|\calE_G - \calE_T|$ is related to the quadrature error
of an integrand, for example, $F(\bsy) := (G(\bsy) -
G_\theta^{[L]}(\bsy))^2$ when $N_{\rm obs} = 1$. Since the quadrature
error for a QMC rule is bounded by the so-called \emph{worst case error}
of the point set, multiplied by the \emph{norm} of the integrand in an
appropriate \emph{weighted} function space, we consider the norm of the
DNN and the norm of the squared difference $F$. In Theorem~\ref{thm:norm}
we provide bounds of the mixed-derivative-based norms for DNNs. In
Theorem~\ref{thm:diff} we prove the key result that, if we \emph{restrict
the network parameters of the DNN to match its derivative bounds to those
of the target function}, then the norm of this squared difference $F$ is
bounded in effectively the same form. Theorem~\ref{thm:err} then shows
that we can construct tailored lattice points to achieve a guaranteed
error bound for the generalization gap, with the implied constant
independent of the dimensionality $s$ of the input data.

Specifically, for a multiindex $\bsnu = (\nu_1,\ldots,\nu_s) \ne\bszero$,
where each $\nu_j$ specifies the number of partial derivatives to be taken
respect to the $j$-th variable of the input data $\bsy$, we find that the
$\bsnu$-th order mixed derivative of the non-periodic DNN
\eqref{eq:DNN-np} in all components of the observables has an upper bound
of the form
\[
  C\,|\bsnu|!\,\bsb^\bsnu =
  C\, \Big(\sum_{j=1}^s \nu_j\Big)!\, \prod_{j=1}^s b_j^{\nu_j},
\]
for some constant $C>0$ and some sequence $(b_j)_{j\ge 1}$ matching the
corresponding derivative bound for a target function. A more complicated
upper bound is obtained for the mixed derivatives of the periodic DNN
\eqref{eq:DNN-per} and a periodic target function. In the latter periodic
case, we obtain the generalization error bound
\begin{align*} 
  \calE_G
  \le {\tt tol} + \calO\big(N^{-r/2}\big),
  \qquad r := \frac{1}{p^*},
\end{align*}
where ${\tt tol}\in (0,1)$ is the tolerance achieved by the training error
in practice, the implied constant in the big-$\calO$ bound is independent
of the input dimensionality~$s$, and $p^* \in (0,1)$ is the
``\emph{summability exponent}'' satisfying $\sum_{j\ge 1} b_j^{p^*} <
\infty$. A smaller~$p^*$ means a faster decay of the sequence $(b_j)_{j\ge 1}$ and a higher
order of convergence.

Section~\ref{sec:pde} provides two example classes of target functions: a
parametric diffusion problem and an algebraic equation mimicking the PDE
problem. In Section~\ref{sec:num} we propose a tailored regularization in
order to ``encourage'' the network parameters to match the desired
derivative bounds as prescribed in Theorem~\ref{thm:diff}. Then as an
early experiment we train DNNs for an algebraic example target function.
Our results indicate that the DNNs trained with tailored regularization
perform significantly better than those trained with only the standard
$\ell_2$ regularization. Section~\ref{sec:proof} contains all the
technical proofs of this paper. The real-valued induction proofs for the
explicit regularity bounds of DNNs are new contributions.
Finally, Section~\ref{sec:conc} provides a conclusion.

\subsection{Related work}

There has been an explosion of works on the theory of DNNs in the last few
years. In this paper we focus on the scenario where the target functions
are smooth, but have very many variables and are costly to evaluate. We
assume that there is a sequence $(b_j)_{j\ge 1}$ which describes the
regularity features of the target function and there is a summability
exponent $p^* \in (0,1)$ for which $\sum_{j\ge 1} b_j^{p^*} < \infty$, as
mentioned above. This is very similar to the setting of the recent review
\cite{ABDM24} which assumes that the target functions are
\emph{holomorphic} (see e.g., \cite{SZ19}): there is a sequence
$(b_j)_{j\ge 1}$ which controls the \emph{anisotropy} of the target
functions, and the results were stated in terms of the summability
exponent $p^*$.

As explained in \cite{ABDM24}, much of the recent research on the
approximation theory of DNNs aim to establish \emph{existence theorems}:
they prove the \emph{expressivity} or \emph{universality} of DNNs, i.e.,
the existence of DNNs of a certain depth and width that approximate a
class of functions to a prescribed accuracy. Such proofs are typically
achieved by \emph{emulation}, i.e., by showing that the network parameters
can be handcrafted so that the DNNs become arbitrarily close to a known
approximation method such as the \emph{best $N$-term approximation}. For
example, for the class of holomorphic functions mentioned above, there
exist DNNs of width $\calO(N^2)$ and depth $\calO(\log(N))$ that achieve
the near-optimal generalization errors of $\calO(N^{-(1/p^*-1/2)})$, see
e.g., \cite[Theorem~7.2]{ABDM24}.

There is however still a significant gap between the approximation theory
of DNNs and their practical performance (see e.g., \cite{AD21,GV24}). Some
papers tried to narrow this gap by developing \emph{practical existence
theory} (see e.g., \cite{ABDM25}) which takes into account the practical
training strategies in minimizing a loss function, including the existence
of a good regularization parameter, see e.g., \cite[Theorem~8.1]{ABDM24}.
For further related literature we refer the readers to the review
\cite{ABDM24} and the references therein.

The review \cite{ABDM24} is devoted to the situation where the learning
methods are independent of the sequence $(b_j)_{j\ge1}$; it is referred to
as the ``unknown anisotropy setting''. In contrast, in this paper we
propose to make explicit use of our knowledge of the sequence
$(b_j)_{j\ge1}$ in our construction of lattice training points as well as
in the design of our practical regularization term. This is a realistic
setting in the examples of parametric PDEs in uncertainty quantification,
since we know explicitly how the random fields are modeled.

``Exploding and vanishing gradients'' in neural network training have been recognized early
on, leading to various approaches to address the training stability issue. Among these, the
most prominent ones are the introduction of gating to enable a constant error flow
during training \cite{LSTM97}, identity mappings \cite{StableNN,ResNet1,ResNet2}
as a simplistic and most natural way of ensuring constant error flow, and numerous normalization
techniques, for example, the currently most popular layer normalization \cite{LayerNorm} to normalize the
inputs to a neuron.

In this paper, we are interested in understanding and developing the basic theory of multi-layer perceptrons (DNNs). Hence, we leave the investigation
of architecture modification for future work and instead turn to the design of the loss function.
Equipped with the regularity bounds in Section~\ref{sec:reg} and error analysis in Section~\ref{sec:err},
we can tailor a loss function that improves multi-layer perceptrons without architectural changes.

\section{Regularity bounds for DNNs} \label{sec:reg}

\subsection{Key definitions, assumptions, and notations}

To establish regularity bounds for both variants of the DNNs
\eqref{eq:DNN-np} and \eqref{eq:DNN-per}, it is convenient to work with a
recursive definition:
\begin{itemize}
\item [\textnormal{(a)}] The non-periodic DNN \eqref{eq:DNN-np} can be
    defined recursively by
\begin{align} \label{eq:rec-np}
  \begin{cases}
  G_\theta^{[0]}(\bsy) := W_0\,\bsy + \bsv_0, \\
  G_\theta^{[\ell]}(\bsy) := W_\ell\, \sigma ( G_\theta^{[\ell-1]}(\bsy)) + \bsv_\ell \quad\mbox{for}\quad \ell\ge 1.
  \end{cases}
\end{align}
\item [\textnormal{(b)}] The periodic DNN \eqref{eq:DNN-per} can be
    defined recursively by
\begin{align} \label{eq:rec-per}
  \begin{cases}
  G_\theta^{[0]}(\bsy) := W_0\,\sin(2\pi\bsy) + \bsv_0, \\
  G_\theta^{[\ell]}(\bsy) := W_\ell\, \sigma ( G_\theta^{[\ell-1]}(\bsy)) + \bsv_\ell
  \quad\mbox{for}\quad \ell\ge 1.
  \end{cases}
\end{align}
\end{itemize}
These recursive definitions split the compositions differently as compared to
how hidden layers are typically defined (layer $\ell$ is given by
$\sigma(W_{\ell-1}\, \bsy + \bsv_{\ell-1})$ for $\bsy\in\bbR^{d_\ell\times
1}$), but \eqref{eq:rec-np} and \eqref{eq:rec-per} are more convenient
for our induction proofs later. 
The view of DNNs as in \eqref{eq:rec-np} and \eqref{eq:rec-per} allows
one to relate the DNN layers to integral operators as in \cite[Sec.~2]{KV23}.

In both cases we will assume that there exist positive sequences
$(\beta_j)_{j\ge 1}$, $(R_\ell)_{\ell\ge 1}$, $(A_n)_{n\ge 1}$ for which
the following three key assumptions hold:
\begin{itemize}
\item The columns of the matrix $W_0$ have bounded vector
    $\infty$-norms
\begin{align} \label{eq:beta}
  \|W_{0,:,j}\|_\infty := \max_{1\le p\le d_1} |W_{0,p,j}| \;\le\; \beta_j
  \quad\mbox{for all}\quad j = 1,\ldots,s.
\end{align}
\item The matrices $W_\ell$ for $\ell\ge 1$ have bounded matrix
    $\infty$-norms (maximum absolute row sum)
\begin{align} \label{eq:R}
  \|W_\ell\|_\infty :=
  \max_{1\le p\le d_{\ell+1}} \sum_{q=1}^{d_\ell} |W_{\ell,p,q}| \;\le\; R_\ell
  \quad\mbox{for all}\quad \ell\ge 1.
\end{align}
\item The activation function $\sigma:\bbR\to\bbR$ is smooth and
    its derivatives satisfy
\begin{align} \label{eq:sigma}
  \|\sigma^{(n)}\|_\infty := \sup_{x\in\bbR} |\sigma^{(n)}(x)| \;\le\; A_n
  \quad\mbox{for all}\quad n\ge 1.
\end{align}
\end{itemize}
The assumptions \eqref{eq:beta}~and~\eqref{eq:R} are analogous to the
assumptions in \cite[Propositions~3.2 and~3.8]{LMRS21}. The assumption
\eqref{eq:sigma} allows for fairly generic smooth activation
functions, including the logistic sigmoid, tanh, and swish (see \eqref{eq:cases}
ahead).

Let $\bbN_0 := \{0,1,2,\ldots\}$ be the set of nonnegative integers. For a
multiindex $\bsnu = (\nu_j)_{j\ge 1} \in \bbN_0^\infty$ we define its
order to be $|\bsnu| := \sum_{j\ge 1} \nu_j$, and we consider the index
set $\indx := \{\bsnu \in \bbN_0^\infty : |\bsnu| < \infty\}$. We write
$\bsm\le\bsnu$ when $m_j\le \nu_j$ for all $j\ge 1$, and
$\binom{\bsnu}{\bsm} := \prod_{j\ge 1} \binom{\nu_j}{m_j}$ using the
notation of binomial coefficients. For any sequence $\bsb = (b_j)_{j\ge
1}$ we write $\bsb^\bsnu := \prod_{j\ge 1} b_j^{\nu_j}$. For any subset
$\setu$ of indices we write $\bsnu_\setu = (\nu_j)_{j\in\setu}$,
$\binom{\bsnu_\setu}{\bsm_\setu} := \prod_{j\in\setu} \binom{\nu_j}{m_j}$,
and $\bsb_\setu^{\bsnu_\setu} := \prod_{j\in\setu} b_j^{\nu_j}$. We define
the mixed partial derivative operator with respect to the parametric
variables $\bsy$ by $\partial^\bsnu := \prod_{j\ge 1}
(\frac{\partial}{\partial y_j})^{\nu_j}$ and $\partial^{\bsnu_\setu} :=
\prod_{j\in\setu} (\frac{\partial}{\partial y_j})^{\nu_j}$.

It turns out that the regularity of the periodic DNN \eqref{eq:rec-per}
depends on the \emph{Stirling number of the second kind}. They are defined
for integers $n$ and $k$ by $\calS(n,0):=\delta_{n,0}$ (i.e., it is 1 if $n=0$
and is $0$ otherwise), $\calS(n,k) := 0$ for $k> n$, and otherwise
\[
  \calS(n,k) := \frac{1}{k!}\sum_{i=0}^k (-1)^{k-i} \binom{k}{i} i^n, \qquad n\ge k.
\]
In combinatorics, $\calS(n,k)$ is the number of ways to partition a set of $n$
objects into $k$ nonempty subsets. In particular, we have $\calS(0,0) = 1$,
$\calS(n,0) = 0 \;\forall n\ge 1$, $\calS(n,1) = 1 \;\forall n\ge 1$, and $\calS(n,n)
= 1 \;\forall n\ge 0$. We introduce the shorthand notation
\[
 \calS(\bsnu,\bsm) := \prod_{j\ge 1} \calS(\nu_j,m_j)
 \quad\mbox{and}\quad
 \calS(\bsnu_\setu,\bsm_\setu) := \prod_{j\in\setu} \calS(\nu_j,m_j).
\]

\subsection{Regularity results}

We are ready to present our regularity results. In the following, we will
use square brackets $[\ell]$ in superscript to indicate all
depth-dependent quantities. Our results are obtained by real-valued
induction argument. As far as we are aware, the explicit regularity bounds
are new.

\begin{theorem} [Regularity bounds for DNNs] \label{thm:der}
Let the sequences $(\beta_j)_{j\ge 1}$, $(R_\ell)_{\ell\ge 1}$,
$(A_n)_{n\ge 1}$ be defined as in \eqref{eq:beta}, \eqref{eq:R},
\eqref{eq:sigma}, respectively. For any depth $\ell\ge 1$, 
any component $1\le p\le d_{\ell+1}$ of $G_\theta^{[\ell]}(\bsy)$,
and any multiindex $\bsnu\in\calI$ with $\bsnu\ne\bszero$, 
we have the following regularity bounds:
\begin{itemize}
\item [\textnormal{(a)}] The non-periodic DNN defined in
    \eqref{eq:rec-np} satisfies
\begin{align} \label{eq:der-np}
  |\partial^\bsnu G_\theta^{[\ell]}(\bsy)_p|
  \le R_\ell\,\bsbeta^{\bsnu}\,\Gamma_{|\bsnu|}^{[\ell]}.
  \qquad\qquad\qquad\qquad\quad\;\,
\end{align}
\item [\textnormal{(b)}] The periodic DNN defined in
    \eqref{eq:rec-per} satisfies
\begin{align} \label{eq:der-per}
  |\partial^\bsnu G_\theta^{[\ell]}(\bsy)_p|
  \le R_\ell\,(2\pi)^{|\bsnu|}
  \sum_{\bsm\le\bsnu} \bsbeta^{\bsm}\,\Gamma_{|\bsm|}^{[\ell]}\,\calS(\bsnu,\bsm).
\end{align}
\end{itemize}
In both cases the sequence $\Gamma_n^{[\ell]}$ is defined recursively by
\begin{align}
  &\begin{cases} \label{eq:G-def}
  \Gamma_n^{[1]} := A_n & \mbox{for $n\ge 1$}, \\
  \Gamma_n^{[\ell]}
  := \displaystyle\sum_{\lambda=1}^n A_\lambda\,R_{\ell-1}^\lambda\,
  \bbB_{n,\lambda}^{[\ell-1]}
  & \mbox{for $\ell\ge 2$ and $n\ge1$,}
  \end{cases} 
  \\
  &\begin{cases} \label{eq:B-def}
  \bbB_{n,1}^{[\ell]} := \Gamma_n^{[\ell]} & \mbox{for $\ell\ge 1$ and $n\ge 1$}, \\
  \bbB_{n,\lambda}^{[\ell]} :=
  \displaystyle\sum_{i=\lambda-1}^{n-1} \binom{n-1}{i}\, \Gamma_{n-i}^{[\ell]}\,\bbB_{i,\lambda-1}^{[\ell]}
  & \mbox{for $\ell\ge 1$ and $n\ge\lambda\ge 2$}. 
  \end{cases}
\end{align}
\end{theorem}

\begin{proof}
See Section~\ref{sec:reg-np} ahead for the proof of \eqref{eq:der-np} and
Section~\ref{sec:reg-per} for the proof of~\eqref{eq:der-per}.
\end{proof}

The regularity bounds \eqref{eq:der-np} and \eqref{eq:der-per} show
structured dependence of the DNNs on the columns of the matrix $W_0$ (via
the sequence $\beta_j$ in~\eqref{eq:beta}) and the subsequent matrices
$W_\ell$ (via the sequence $R_\ell$ in~\eqref{eq:R}) as well as on the
activation function $\sigma$ (via the sequence $A_n$ in~\eqref{eq:sigma}).
In particular, the sequence $\Gamma_n^{[\ell]}$ is determined by the
combined effect of the activation function $\sigma$ (via the sequence
$A_n$) and the interior matrices $W_1,\ldots,W_{\ell-1}$ (via
the numbers $R_1,\ldots,R_{\ell-1}$ in~\eqref{eq:R}).

When the dimensionality of the input parameters $d_0 = s$ is large, it
may make sense to impose a condition that the sequence $\beta_j$ is
ordered and decays sufficiently fast, i.e., the columns of the matrix
$W_0$ should decay in vector $\infty$-norms. This is consistent with the
accepted wisdom in contemporary high-dimensional QMC integration theory
that the underlying function space setting should be ``weighted'' so that
the parameters or the integration variables have decaying importance.

Considering mixed derivatives with respect to the input parameters $\bsy$
with high order $|\bsnu|$, it makes sense to
restrict the matrix elements of $W_\ell$ such that their matrix
$\infty$-norms are bounded by a reasonably small constant.
At the same time, we would also want to make sure that the derivatives of
the activation function $\sigma$ do not grow too quickly in magnitude so
that the growth of the sequence $\Gamma_n^{[\ell]}$ can be sufficiently
compensated by the decay of the sequence $\beta_j$. 

In Section~\ref{sec:act} we consider the following activation
functions and determine explicit expressions for $A_n$ in
\eqref{eq:sigma}:
\begin{align} \label{eq:cases}
\begin{cases}
\begin{tabular}{lll}
 sigmoid: & $\sigma(x) = \displaystyle\frac{1}{1+e^{-x}}$,
 & $A_n = n!$\,, 
 \vspace{0.2cm} \\
 tanh: & $\sigma(x) = \displaystyle\frac{e^x-e^{-x}}{e^x+e^{-x}}$,
 & $A_n = 2^n\,n!$\,, 
 \vspace{0.2cm} \\
 swish: & $\sigma(x) = \displaystyle\frac{x}{1+e^{-x}}$,
 & $A_n = 1.1\,n!$\,.
\end{tabular}
\end{cases}
\end{align}
The ReLU activation function $\sigma(x) = \max\{x,0\}$ is not smooth and so condition
\eqref{eq:sigma} does not hold.
However, a generalization of swish, namely
$x/(1+e^{-c x})$ for $c \rightarrow \infty$ converges to the
ReLU activation function.  Like ReLU, the swish function is unbounded as $x\to+\infty$.

For the remainder of this paper, we will illustrate our results for
activation functions whose derivatives are bounded in the sense of \eqref{eq:sigma} by
the common form
\begin{align} \label{eq:common}
  A_n = \xi\,\tau^n\,n!
  \qquad\mbox{for some $\xi>0$ and $\tau>0$}.
\end{align}
Clearly the sigmoid case corresponds to $\xi = \tau = 1$, the tanh case
corresponds to $\xi=1$ and $\tau=2$, and the swish case corresponds to
$\xi = 1.1$ and $\tau=1$.

\begin{theorem}[Regularity bounds for DNNs with $A_n$ of the form \eqref{eq:common}] \label{thm:reg}
Let the sequences $(\beta_j)_{j\ge 1}$, $(R_\ell)_{\ell\ge 1}$,
$(A_n)_{n\ge 1}$ be defined as in \eqref{eq:beta}, \eqref{eq:R},
\eqref{eq:sigma}, respectively, with $A_n$ given by \eqref{eq:common}. For
any depth $\ell\ge 1$, any component $1\le p\le d_{\ell+1}$, and any
multiindex $\bsnu\in\calI$ $($including $\bsnu=\bszero$$)$, we have the
following regularity bounds:
\begin{itemize}
\item [\textnormal{(a)}] The non-periodic DNN defined in
    \eqref{eq:rec-np} satisfies
\begin{align} \label{eq:reg-np}
  |\partial^\bsnu G_\theta^{[\ell]}(\bsy)_p|
  \le C_\ell\, \,|\bsnu|!\, (S_\ell\,\bsbeta)^{\bsnu}.
  \qquad\qquad\qquad\qquad\quad
\end{align}

\item [\textnormal{(b)}] The periodic DNN defined in
    \eqref{eq:rec-per} satisfies
\begin{align} \label{eq:reg-per}
  |\partial^\bsnu G_\theta^{[\ell]}(\bsy)_p|
  \le C_\ell\, (2\pi)^{|\bsnu|}
  \sum_{\bsm\le\bsnu} |\bsm|!\,(S_\ell\,\bsbeta)^{\bsm}\,\calS(\bsnu,\bsm).
\end{align}
\end{itemize}
In both cases we define 
\begin{align} \label{eq:CSP}
  &C_\ell := \max\bigg\{\displaystyle
 \|G^{[\ell]}_\theta\|_{\infty},\;
 \frac{P_\ell}{S_\ell} \bigg\},
  \quad
  S_\ell := \tau \sum_{k=0}^{\ell-1} P_k, 
  \\
  &P_0 := 1, \quad
  P_k := \displaystyle\prod_{t=1}^k (\xi\,\tau\, R_t)
  \quad\mbox{for $k\ge 1$}. \nonumber
\end{align}
\end{theorem}

\begin{proof}
See Lemma~\ref{lem:complex} in Section~\ref{sec:act} for a proof that
if $A_n = \xi\,\tau^n\,n!$ then the sequence $\Gamma_n^{[\ell]}$ defined
recursively in \eqref{eq:G-def}--\eqref{eq:B-def} is given explicitly for $\ell\ge 1$ by
\begin{align*}
 \Gamma_n^{[\ell]}
 &= P_{\ell-1}\,\bigg(\sum_{k=0}^{\ell-1} P_k\bigg)^{n-1}\,\xi\,\tau^n\,n!
 = (\xi\,\tau\,P_{\ell-1})\,\bigg(\tau \sum_{k=0}^{\ell-1} P_k\bigg)^{n-1}\,n!
 = \frac{\xi\,\tau\,P_{\ell-1}}{S_\ell}\,S_\ell^n\,n!\,.
\end{align*}
Substituting this into Theorem~\ref{thm:der} and then using
$R_\ell\,(\xi\,\tau\,P_{\ell-1})/S_\ell = P_\ell/S_\ell \le C_\ell$ and
$S_\ell^{|\bsm|}\, \bsbeta^\bsm = (S_\ell\,\bsbeta)^\bsm$, we obtain the
bounds \eqref{eq:reg-np} and \eqref{eq:reg-per} for $\bsnu\ne\bszero$. The
case $\bsnu=\bszero$ holds trivially since $|G_\theta^{[\ell]}(\bsy)_p|
\le \|G_\theta^{[\ell]}\|_\infty \le C_\ell$.
\end{proof}

It is interesting to see that the regularity bounds \eqref{eq:reg-np}
and \eqref{eq:reg-per} are of the POD (``product and order dependent'')
and SPOD (``smoothness-driven product and order dependent'') forms which
we have encountered in previous works on parametric PDEs in uncertainty
quantification, see \cite{KSS12} and \cite{DKLNS14} where they first
appeared, as well as, e.g., \cite{HHKKS24,KKKNS22,KKS20,KKS24}.

Furthermore, we observe from the regularity bounds \eqref{eq:reg-np} and
\eqref{eq:reg-per} that, as the depth $\ell$ increases, we would need to
restrict the magnitude of $R_k$ for $1\le k\le \ell-1$ in order for the sum $S_\ell$ to be bounded.

\section{Error analysis} \label{sec:err}

\subsection{Preliminaries}

The goal of the training process in machine learning is to find the
network parameters $\theta\in\Theta$ that minimize a \emph{loss
function}, in our case
\[
  \calJ :=
  \calJ(\theta) := \frac{1}{N} \sum_{k=1}^N
  \|G(\bsy_k) - G_{\theta}^{[L]}(\bsy_k)\|_2^2,
\]
where we took the vector $2$-norm with respect to the $N_{\rm obs}$ components in the observables.
It is common to regularize the minimization problem and so we use
\begin{align} \label{eq:minimization}
  \theta^* := \argmin_{\theta\in\Theta}\,
  \big(\calJ(\theta) + \lambda\, \calR(\theta) \big),
\end{align}
with an appropriately chosen regularization term $\calR(\theta)$.

The error of interest is the so-called \emph{generalization error} (or
population risk) which measures the error of the network on unseen data.
This is defined as
\begin{align} \label{eq:EG}
  \calE_G := \calE_G(\theta) := \| G - G_{\theta}^{[L]} \|_{L_2(Y)}
  = \bigg(\int_Y \| G(\bsy) - G_{\theta}^{[L]}(\bsy) \|_2^2 \,\rd\bsy\bigg)^{1/2}.
\end{align}
There is also the computable \emph{training error} (or empirical risk),
given by
\begin{align} \label{eq:ET}
  \calE_T := \calE_T(\theta) :=
  \bigg(\frac{1}{N} \sum_{k=1}^N \|G(\bsy_k) - G_{\theta}^{[L]}(\bsy_k) \|_2^2 \bigg)^{1/2},
\end{align}
leading to the so-called \emph{generalization gap} $|\calE_G - \calE_T|$.
Obviously, for any $\calE_G, \calE_T \ge 0$ we have
\begin{align} \label{eq:gap}
  \calE_G \le \calE_T + |\calE_G - \calE_T|.
\end{align}
If $\calE_G \le\calE_T$ then the generalization error is already under control. On the other hand, if $\calE_G >\calE_T$ then we have 
$\max(\calE_T,\calE_G - \calE_T) \le \calE_G = \calE_T + (\calE_G - \calE_T)$,
and it is the dominant one between $\calE_T$ and $(\calE_G - \calE_T)$
that will capture the trend of $\calE_G$. This explains why we want to control the gap $|\calE_G - \calE_T|$ relative to the training error $\calE_T$.

The square roots in \eqref{eq:EG} and \eqref{eq:ET} make it difficult to analyze the generalization gap $|\calE_G - \calE_T|$ directly. However, we observe that the difference $\calE_G^2 - \calE_T^2$ is the
quadrature error for the integrand $F(\bsy) = \|(G(\bsy) -
G_{\theta}^{[L]}(\bsy)\|_2^2$. For example, when $N_{\rm obs} = 1$, the integrand is simply
$F(\bsy) = (G(\bsy) - G_{\theta}^{[L]}(\bsy))^2$.

The paper \cite{MR21} considered $N_{\rm obs}=1$ and the $1$-norm
instead of the $2$-norm in the definitions of \eqref{eq:EG} and \eqref{eq:ET}. Interpreting the generalization gap as a quadrature error for the integrand $F(\bsy) = |G(\bsy) -
G_{\theta}^{[L]}(\bsy)|$, they used a smooth function to replace the absolute
value function in the neighborhood of the origin, together with the
Koksma--Hlawka inequality for low-discrepancy sequences, to conclude that
$|\calE_G - \calE_T| \le c\,(\log N)^s/N + 2\delta$, where $\delta>0$
is their tolerance parameter used in mollifying the absolute value
function.

The paper \cite{LMRS21} considered the $2$-norm together with
higher order QMC rules to obtain
\begin{align} \label{eq:err1}
  |\calE_G - \calE_T| = \frac{|\calE_G^2 - \calE_T^2|}{\calE_G + \calE_T}
  \le \frac{|\calE_G^2 - \calE_T^2|}{\calE_G}
  \le \frac{|\calE_G^2 - \calE_T^2|}{\inf_{\theta\in\Theta} \calE_G(\theta)}
  \le c_{\rm arch}\, N^{-r}, 
\end{align}
where the finiteness of $c_{\rm arch}$ relies on an additional assumption,
see \cite[formula~(3.9)]{LMRS21}, which depends on the network
architecture. An easy alternative bound in \cite{LMRS21} is given by
\begin{align} \label{eq:err2}
  |\calE_G - \calE_T| \le \sqrt{|\calE_G^2 - \calE_T^2|}
  \le c\, N^{-r/2} 
\end{align}
with a constant $c$ independent of the network architecture, but at only
half the convergence rate. In both cases, $r = 1/p^*$ where $p^*$ is the
summability exponent of a sequence $(b_j)_{j\ge 1}$ characterising the
holomorphic target functions, as mentioned in the introduction.

In this paper we will instead use lattice rules in both the non-periodic
and periodic settings. Furthermore, instead of the holomorphic argument,
we will use the real-valued regularity bounds from the previous section.

\subsection{Our error bound}

In general, the \emph{worst case error} for any QMC rule $Q(F) = \frac{1}{N} \sum_{k=1}^N F(\bst_k)$ approximating the integral $I(F)=\int_Y F(\bsy)\,\rd\bsy$ of all functions $F$ in a normed space $\calW$ is defined by
\[
  e^{\rm wor}(Q,\calW) := \sup_{F\in\calW,\,\|F\|_{\calW}\le 1} | I(F) - Q(F) |.
\]
Since $Q$ is linear, it follows that $|I(F) - Q(F)| \le e^{\rm wor}(Q,\calW)\, \|F\|_\calW$ for all $F\in\calW$, i.e., the quadrature error is bounded by the product of the worst case error and the norm of the integrand. The worst case error depends on the QMC rule $Q$ as well as properties of the function space $\calW$.
In what follows we pick a lattice rule~\eqref{eq:lat1} for the QMC rule which is fully determined by its generating vector $\bsz$ and the number of points $N$.
We will hence denote its worst case error in a space $\calW$ by $e^{\rm wor}_N(\bsz, \calW)$.

Specifically, for an arbitrary number of observables $N_{\rm obs} = d_{L+1}$, our error
expression of interest (see \eqref{eq:err1} and \eqref{eq:err2}) satisfies
\begin{align} \label{eq:wce}
  |\calE_G^2 - \calE_T^2|
  &= \bigg| \int_Y \big\| G(\bsy) - G_{\theta}^{[L]}(\bsy) \big\|_2^2 \,\rd\bsy
  - \frac{1}{N} \sum_{k=1}^N \big\|G(\bsy_k) - G_{\theta}^{[L]}(\bsy_k)\big\|_2^2\, \bigg| \nonumber\\
  &= \bigg|
  \int_Y \sum_{p=1}^{N_{\rm obs}} \big( G(\bsy)_p - G_{\theta}^{[L]}(\bsy)_p \big)^2 \,\rd\bsy
  - \frac{1}{N} \sum_{k=1}^N \sum_{p=1}^{N_{\rm obs}}
  \big(G(\bsy_k)_p - G_{\theta}^{[L]}(\bsy_k)_p\big)^2
  \bigg| \nonumber\\
  &\le \sum_{p=1}^{N_{\rm obs}} \bigg|
  \int_Y \big( G(\bsy)_p - G_{\theta}^{[L]}(\bsy)_p \big)^2 \,\rd\bsy
  - \frac{1}{N} \sum_{k=1}^N \big(G(\bsy_k)_p - G_{\theta}^{[L]}(\bsy_k)_p\big)^2
  \bigg| \nonumber\\
  &\le \sum_{p=1}^{N_{\rm obs}} e^{\rm wor}_N(\bsz, \calW_{\alpha,\bsgamma})\,
  \big\| \big(G(\cdot)_p - G_{\theta}^{[L]}(\cdot)_p\big)^2 \big\|_{\calW_{\alpha,\bsgamma}},
\end{align}
where for each component $p$ we have the \emph{worst case integration
error} $e^{\rm wor}_N(\bsz, \calW_{\alpha,\bsgamma})$ for a lattice rule (with $N$ points and a
generating vector $\bsz$) in a \emph{weighted function space}
$\calW_{\alpha,\bsgamma}$, multiplied by the norm of the scalar-valued
integrand $F(\bsy) := (G(\bsy)_p - G_{\theta}^{[L]}(\bsy)_p)^2$ in this
function space. This space has an integer \emph{smoothness parameter}
$\alpha$ corresponding to the order of derivatives available in each
variable, together with positive \emph{weight parameters} $\bsgamma =
(\gamma_\setu)_{\setu\subseteq\{1:s\}}$, $\gamma_\emptyset:=1$, which are
used to moderate the relative importance between subsets of variables,
see three settings in the next subsection. The worst case error depends on the generating vector $\bsz$ as well as the smoothness parameter $\alpha$ and the weights $\bsgamma$.
To apply the theory for
lattice rules, we need to show that each integrand $(G(\cdot)_p -
G_{\theta}^{[L]}(\cdot)_p)^2$ indeed belongs to the appropriate weighted
function space $\calW_{\alpha,\bsgamma}$ and has a nice bound on its norm.

Note that after the minimization \eqref{eq:minimization} the network
parameters $\theta \approx \theta^*$ depend on $N$ and the training
points, as well as other hyperparameters in the network architecture. In
summary, putting \eqref{eq:gap}, \eqref{eq:err2}, \eqref{eq:wce} together,
the generalization error for the trained DNN satisfies
\begin{align} \label{eq:combine}
  \calE_G \le \calE_T + |\calE_G - \calE_T|
  &\le \calE_T + \sqrt{|\calE_G^2 - \calE_T^2|} \nonumber\\
  &\le \calE_T +
  \sqrt{
  \textstyle\sum_{p=1}^{N_{\rm obs}} e^{\rm wor}_N(\bsz, \calW_{\alpha,\bsgamma})\,
  \| (G(\cdot)_p - G_{\theta}^{[L]}(\cdot)_p)^2 \|_{\calW_{\alpha,\bsgamma}}
  }.
\end{align}
Hence, assuming that the training error $\calE_T$ is under control (i.e.,
small relative to the generalization gap), we have a guaranteed
convergence with respect to $N$ from a good lattice rule, as long as the
norms $\| (G(\cdot)_p - G_{\theta}^{[L]}(\cdot)_p)^2
\|_{\calW_{\alpha,\bsgamma}}$ are bounded. As $N$ increases we expect
$G_{\theta}^{[L]}$ to become closer to $G$ in some pointwise sense, and these norms in
\eqref{eq:combine} could also become smaller with increasing~$N$, but we
do not know how to exploit this. This situation is different from a
typical QMC integration problem where the integrand is fixed and
independent of $N$.

\subsection{Review of three theoretical settings for lattice rules} \label{sec:qmc}

There are known algorithms (``component-by-component construction'') for
finding good lattice generating vectors $\bsz$ that achieve nearly the
optimal convergence rates for the worst case errors $e^{\rm wor}_N(\bsz, \calW_{\alpha,\bsgamma})$
in \eqref{eq:wce} in appropriately paired function space settings.
Sometimes we consider ``randomly shifted'' lattice rules and then the
error bounds hold in a root-mean-square sense with respect to the random
shift. We illustrate three possible theoretical settings in this
paper, see e.g.,~\cite[Theorem~5.8]{DKS13} for (a), \cite[Theorem~5.12]{DKS13} for (b), and \cite[Lemma~3.1]{KKS20} for (c).
(Note that there are lots of other theoretical results for lattice rules and other QMC rules.)

\begin{itemize}
\item [\textnormal{(a)}] For a ``randomly shifted'' lattice rule in a
    (non-periodic, ``unanchored'') weighted Sobolev space
    ${\calW_{1,\bsgamma}}$ with smoothness parameter $\alpha = 1$, we
    set $\bsnu_\setu = (1,\ldots,1)$ to define
\begin{align} \label{eq:norm-np}
 \|F\|_{\calW_{1,\bsgamma}}^2
 := \sum_{\setu\subseteq\{1:s\}} \frac{1}{\gamma_\setu}
 \int_{[0,1]^{|\setu|}} \bigg| \int_{[0,1]^{s-|\setu|}}
 \partial^{\bsnu_\setu} F(\bsy)\,
 \rd\bsy_{\{1:s\}\setminus\setu}\bigg|^2\rd\bsy_\setu.
\end{align}
A good generating vector $\bsz$ can be constructed such that the
root-mean-square worst case error satisfies
\begin{align} \label{eq:wce-np}
 {\rm r.m.s.}\; e^{\rm wor}_N(\bsz, \calW_{1,\bsgamma})
  \le \bigg(
  \frac{2}{N} \sum_{\emptyset\ne \setu \subseteq \{1:s\}} \gamma_\setu^\lambda\,
  \bigg[
  \frac{2\zeta(2\lambda)\,2^\lambda}{(2\pi)^{2\lambda}} 
  \bigg]^{|\setu|}\bigg)^{\frac{1}{2\lambda}}\,
  \quad\forall\;\lambda\in (\tfrac{1}{2},1],
\end{align}
where $\zeta(x) := \sum_{h=1}^\infty h^{-x}$ is the Riemann zeta
function.

\item [\textnormal{(b)}] For a lattice rule in a (periodic, Hilbert)
    weighted Korobov space $\calW_{\alpha,\bsgamma}$ with integer
    smoothness parameter $\alpha \ge 1$, we set $\bsnu_\setu =
    (\alpha,\ldots,\alpha)$ to define
\begin{align} \label{eq:norm-per}
 \|F\|_{\calW_{\alpha,\bsgamma}}^2
 := \sum_{\setu\subseteq\{1:s\}} 
 \frac{1}{\gamma_\setu}
 \int_{[0,1]^{|\setu|}} \bigg| \int_{[0,1]^{s-|\setu|}}
 \partial^{\bsnu_\setu} F(\bsy)\,
 \rd\bsy_{\{1:s\}\setminus\setu}\bigg|^2\rd\bsy_\setu.
\end{align}
A good generating vector $\bsz$ can be constructed such that the worst
case error satisfies
\begin{align} \label{eq:wce-per}
  e^{\rm wor}_N(\bsz, \calW_{\alpha,\bsgamma})
  \le \bigg(
  \frac{2}{N} \sum_{\emptyset\ne \setu\subseteq \{1:s\}}
  \gamma_\setu^\lambda\,
  \bigg[\frac{2\zeta(2\alpha\lambda)}{(2\pi)^{2\alpha\lambda}}
  \bigg]^{|\setu|}
  \bigg)^{\frac{1}{2\lambda}}\, \quad\forall\;\lambda\in
  (\tfrac{1}{2\alpha},1].
\end{align}

\item [\textnormal{(c)}] For a lattice rule in a (periodic,
    non-Hilbert) weighted Korobov space $\calK_{\alpha,\bsgamma}$ with
    integer smoothness parameter $\alpha \ge 2$, we set $\bsnu_\setu =
    (\alpha,\ldots,\alpha)$ and the norm (typically defined in terms
    of Fourier series) is bounded by
\begin{align} \label{eq:norm-K}
 \|F\|_{\calK_{\alpha,\bsgamma}}
 \le \max_{\setu\subseteq\{1:s\}} 
 \frac{1}{\gamma_\setu}
 \int_{[0,1]^{|\setu|}} \bigg| \int_{[0,1]^{s-|\setu|}}
 \partial^{\bsnu_\setu} F(\bsy)\,
 \rd\bsy_{\{1:s\}\setminus\setu}\bigg| \,\rd\bsy_\setu.
\end{align}
A good generating vector $\bsz$ can be constructed such that the worst
case error satisfies
\begin{align} \label{eq:wce-K}
  e^{\rm wor}_N(\bsz, \calK_{\alpha,\bsgamma})
  \le \bigg(
  \frac{2}{N} \sum_{\emptyset\ne \setu \subseteq \{1:s\}}
  \gamma_\setu^\lambda\,
  \bigg[\frac{2\zeta(\alpha\lambda)}{(2\pi)^{\alpha\lambda}}
  \bigg]^{|\setu|}
  \bigg)^{\frac{1}{\lambda}}\, \quad\forall\;\lambda\in
  (\tfrac{1}{\alpha},1].
\end{align}
\end{itemize}
All three results above were stated for $N$ being a power of a prime, and
similar results hold for any composite $N$. We can also construct
``embedded'' lattice rules~\cite{CKN06} to be used for a range of values
of~$N$ (e.g., from $2^{10}$ to $2^{20}$), at the price of an extra small
constant factor in the error bound.

We have deliberately stated the above results so that the norm
\eqref{eq:norm-np} for the non-periodic setting (a) matches the norm
\eqref{eq:norm-per} in the periodic setting (b) with $\alpha=1$. The
Sobolev norm with higher smoothness for non-periodic functions will
involve all mixed derivatives up to order~$\alpha$, unlike
\eqref{eq:norm-per} which contains only the highest order derivative. The
Hilbert Korobov norm is typically defined in terms of the decay of Fourier
series with a real parameter $\alpha>1/2$, and only when $\alpha$ is an
integer can the norm be written in the form \eqref{eq:norm-per} involving
mixed derivatives. Sometimes the scaling of $\alpha/2$ is used instead of
$\alpha$, and often the scaling of $(2\pi)^{2\alpha|\setu|}$ is absorbed
into the weights $\gamma_\setu$. Similarly, the non-Hilbert Korobov norm
is typically defined in terms of the decay of Fourier series with a real
parameter $\alpha>1$, and only when $\alpha$ is an integer can the norm
satisfy the bound \eqref{eq:norm-K} involving mixed derivatives. The
corresponding worst case error is known as ``$P_\alpha$'' in the classical
theory for lattice rules \cite{SJ94}. We have denoted this space
differently as $\calK_{\alpha,\bsgamma}$ so that there is no confusion
with the norms \eqref{eq:norm-np} and \eqref{eq:norm-per}.

Taking $\lambda$ arbitrarily close to $1/2$ in \eqref{eq:wce-np}, we see that
setting~(a) can achieve close to $\calO(N^{-1})$ convergence for the root-mean-square worst case error. Similarly, taking $\lambda$ arbitrarily close to $1/(2\alpha)$ and $1/\alpha$ in \eqref{eq:wce-per} and \eqref{eq:wce-K}, respectively, we see that both settings (b) and (c) can achieve close to $\calO(N^{-\alpha})$
convergence for the worst case error, however, the requirements on the weights $\gamma_\setu$ to
ensure that the implied constants are independent of the input
dimensionality $s$ can be quite different. Later in Theorem~\ref{thm:err}, we will see pros and cons of the non-Hilbert setting.

One may notice that the bound \eqref{eq:wce-np} has an extra factor of
$2^\lambda$ inside the square brackets compared to \eqref{eq:wce-per} with
$\alpha=1$; this arose from the technique of averaging with respect to the
random shift and working with an associated ``shift-invariant'' setting.
We can also apply random shifting in the periodic setting and the
resulting root-mean-square error bound would remain the same as
\eqref{eq:wce-per}. Random shifting is useful as it provides a practical
error estimate.

The computational cost for the construction of good lattice generating vectors depends on the structure of weights $\bsgamma$, ranging from $\calO(s\,N\log N)$ for product weights, to $\calO(s\,N\log N + s^2\,N)$ for POD weights, to $\calO(s\,N\log N + \alpha^2\,s^2\,N)$ for SPOD weights, see e.g., \cite{DKS13,KKS20,NC06,N14}. This will be an offline precomputation and the resulting generating vector $\bsz$ can be stored for use later.

\subsection{Norm bounds for DNNs} \label{sec:normbounds}

We see from \eqref{eq:wce} and \eqref{eq:combine} that we need to bound
the norms $\| (G(\cdot)_p - G_{\theta}^{[L]}(\cdot)_p)^2
\|_{\calW_{\alpha,\bsgamma}}$. We begin by obtaining bounds on the norm of
the DNN.

\begin{theorem}[Norm bounds for DNNs] \label{thm:norm}
Let the sequences $(\beta_j)_{j\ge 1}$, $(R_\ell)_{\ell\ge 1}$,
$(A_n)_{n\ge 1}$ be defined as in \eqref{eq:beta}, \eqref{eq:R},
\eqref{eq:sigma}, respectively, with $A_n$ given by \eqref{eq:common}.
For depth $L\ge 1$, let $C_L$ and $S_L$ be
defined as in \eqref{eq:CSP}, and $d_{L+1} = N_{\rm obs}$.

\begin{itemize}
\item [\textnormal{(a)}] The norm \eqref{eq:norm-np} of every
    component $1\le p\le N_{\rm obs}$ of the non-periodic DNN
    \eqref{eq:rec-np} satisfies
\begin{align} \label{eq:GL-np}
 \big\|G^{[L]}_\theta(\cdot)_p\big\|_{\calW_{1,\bsgamma}}^2
 &\le
 C_L^2 \sum_{\setu\subseteq\{1:s\}} \frac{1}{\gamma_\setu} \bigg(
 |\setu|! \prod_{j\in\setu} (S_L\,\beta_j) \bigg)^2.
\end{align}

\item [\textnormal{(b)}] The norm \eqref{eq:norm-per} of every
    component $1\le p\le N_{\rm obs}$ of the periodic DNN
    \eqref{eq:rec-per} satisfies, with $\bsalpha_\setu =
    (\alpha,\ldots,\alpha)$,
\begin{align} \label{eq:GL-per}
 &\big\|G^{[L]}_\theta(\cdot)_p \big\|_{\calW_{\alpha,\bsgamma}}^2 \\
 &\le C_L^2
 \sum_{\setu\subseteq\{1:s\}} \frac{(2\pi)^{2\alpha|\setu|}}{\gamma_\setu} \bigg(
 \sum_{\bsm_\setu\le\bsalpha_\setu}
 |\bsm_\setu|!\, \calS(\bsalpha_\setu,\bsm_\setu)\,
 \prod_{j\in\setu} (S_L\,\beta_j)^{m_j} \bigg)^2. \nonumber
\end{align}

\item [\textnormal{(c)}] The norm \eqref{eq:norm-K} of every component
    $1\le p\le N_{\rm obs}$ of the periodic DNN \eqref{eq:rec-per}
    satisfies, with $\bsalpha_\setu = (\alpha,\ldots,\alpha)$,
\begin{align} \label{eq:GL-K}
 \big\|G^{[L]}_\theta(\cdot)_p \big\|_{\calK_{\alpha,\bsgamma}}
 &\le C_L
 \max_{\setu\subseteq\{1:s\}} \frac{(2\pi)^{\alpha|\setu|}}{\gamma_\setu}\!\!
 \sum_{\bsm_\setu\le\bsalpha_\setu}\!
 |\bsm_\setu|!\, \calS(\bsalpha_\setu,\bsm_\setu)\,
 \prod_{j\in\setu} (S_L\,\beta_j)^{m_j} .  
\end{align}
\end{itemize}
\end{theorem}

\begin{proof}
The three results follow by applying the bounds \eqref{eq:reg-np} and
\eqref{eq:reg-per} in the definition of the norms \eqref{eq:norm-np},
\eqref{eq:norm-per} and \eqref{eq:norm-K}.
\end{proof}

For bounded activation functions (e.g., $\|\sigma\|_{\infty} = 1$ for
sigmoid and tanh), we have $\|G^{[L]}_\theta\|_{\infty} \le
R_L\,\|\sigma\|_{\infty} + \|\bsv_L\|_\infty$, where $\bsv_L$ is the final
vector in \eqref{eq:DNN-np} or \eqref{eq:DNN-per}. Also, we have
\[
 \frac{P_L}{S_L}
 = \frac{(\xi\,\tau\,R_L)\, P_{L-1}}{\tau\sum_{k=0}^{L-1} P_k}
 \le \xi\,R_L.
\]
Thus $C_L$ defined in \eqref{eq:CSP} is bounded by
\[
 C_L\le \max\{R_L\,\|\sigma\|_\infty + \|\bsv_L\|_\infty,\; \xi\,R_L\}.
\]
For unbounded activation functions (e.g., swish) the DNN can be unbounded,
and it will be necessary in practice to impose a stronger restriction on
$R_L$.

We can restrict the network parameters of the DNNs to match the regularity
features of the target functions. In turn this allows us to obtain a
theoretical bound on the norms in \eqref{eq:combine}.

\begin{theorem}[Match regularity of DNNs to observables] \label{thm:diff}
Let the sequences $(\beta_j)_{j\ge 1}$, $(R_\ell)_{\ell\ge 1}$,
$(A_n)_{n\ge 1}$ be defined as in \eqref{eq:beta}, \eqref{eq:R},
\eqref{eq:sigma}, respectively, with $A_n$ given by \eqref{eq:common}.
For depth $L\ge 1$, let $C_L$ and $S_L$ be
defined as in \eqref{eq:CSP}, and $d_{L+1} = N_{\rm obs}$.
Given a sequence $(b_j)_{j\ge 1}$ and a constant $C>0$, we 
\begin{itemize}
\item restrict the elements of the matrices $W_1,\ldots,W_{L-1}$ so
    that for some constant $\rho>0$,
\begin{align} \label{eq:demand1}
 & R_\ell \le \rho
  \quad\mbox{for all $\ell= 1, \ldots,L-1$},
 \\
 \label{eq:demand1b}
 \Big( \mbox{thus}\quad S_L \le \tau\,L \quad\mbox{if}
 &\quad \xi\,\tau\,\rho = 1, \quad\mbox{and}\quad
 S_L \le \tau\frac{(\xi\,\tau\,\rho)^L-1}{\xi\,\tau\,\rho-1}
 \quad\mbox{otherwise} \Big),
\end{align}
\item restrict the elements of the matrix $W_0$ so that
\begin{align} \label{eq:demand2}
  \beta_j \le \frac{b_j}{S_L}
  \quad\mbox{for all $j=1,\ldots,s$},
\end{align}
\item restrict the matrices $W_1,\ldots,W_L$ and vector $\bsv_L$ so that
\begin{align} \label{eq:demand3}
 C_L \le C.
\end{align}
\end{itemize}
\begin{itemize}
\item [\textnormal{(a)}] Suppose a non-periodic target function
    $G(\bsy)$ satisfies the regularity bound: for all multiindices
    $\bsnu\in\calI$ and all components $1\le p\le N_{\rm obs}$,
\begin{align} \label{eq:tar-np}
    |\partial^\bsnu G(\bsy)_p| \le C\,|\bsnu|!\, \bsb^{\bsnu}.
\end{align}
Then the non-periodic DNN \eqref{eq:DNN-np} with restrictions
\eqref{eq:demand1}--\eqref{eq:demand3} satisfies the same regularity
bound \eqref{eq:tar-np}, and the norm \eqref{eq:norm-np} of both
functions satisfy for all components $1\le p\le N_{\rm obs}$,
\begin{align}
 &\max\big\{
 \big\|G(\cdot)_p\big\|_{\calW_{1,\bsgamma}}^2, 
 \big\|G^{[L]}_\theta(\cdot)_p\big\|_{\calW_{1,\bsgamma}}^2 \hphantom{,}
 \big\}
 \le C^2
 \sum_{\setu\subseteq\{1:s\}} \frac{1}{\gamma_\setu} \bigg(
 |\setu|!\, \prod_{j\in\setu} b_j \bigg)^2,  \label{eq:G-np}
 \\
 &\big\|\big(G(\cdot)_p- G^{[L]}_\theta(\cdot)_p\big)^2\big\|_{\calW_{1,\bsgamma}}^2
 \le 16\,C^4
 \sum_{\setu\subseteq\{1:s\}} \frac{1}{\gamma_\setu} \bigg(
 (|\setu|+1)!\, \prod_{j\in\setu} b_j
 \bigg)^2. \label{eq:final-np}
\end{align}

\item [\textnormal{(b)}] Suppose a periodic target function $G(\bsy)$
    satisfies the regularity bound: for all mutiindices
    $\bsnu\in\calI$ and all components $1\le p\le N_{\rm obs}$,
\begin{align} \label{eq:tar-per}
  |\partial^\bsnu G(\bsy)_p|
  \le C\,(2\pi)^{|\bsnu|} \sum_{\bsm\le\bsnu} |\bsm|!\,\bsb^{\bsm}\,\calS(\bsnu,\bsm).
\end{align}
Then the periodic DNN \eqref{eq:DNN-per} with restrictions
\eqref{eq:demand1}--\eqref{eq:demand3} satisfies the same regularity
bound \eqref{eq:tar-per}, and the norm \eqref{eq:norm-per} of both
functions satisfy for all components $1\le p\le N_{\rm obs}$,
\begin{align}
 &\max\big\{
 \big\|G(\cdot)_p\big\|_{\calW_{\alpha,\bsgamma}}^2, 
 \big\|G^{[L]}_\theta(\cdot)_p\big\|_{\calW_{\alpha,\bsgamma}}^2 \hphantom{,}
 \big\} \nonumber\\
 &\le C^2\!\!
 \sum_{\setu\subseteq\{1:s\}} \!\!\!\!
 \frac{(2\pi)^{2\alpha|\setu|}}{\gamma_\setu} \bigg(
 \sum_{\bsm_\setu\le\bsalpha_\setu} |\bsm_\setu|!\,
 \calS(\bsalpha_\setu,\bsm_\setu)\, \prod_{j\in\setu} b_j^{m_j}\!
 \bigg)^2,
 \label{eq:G-per}
 \\
 &\big\|\big(G(\cdot)_p - G^{[L]}_\theta(\cdot)_p\big)^2\big\|_{\calW_{\alpha,\bsgamma}}^2 \nonumber \\
 &\le 16\,C^4\!\!
 \sum_{\setu\subseteq\{1:s\}} \!\!\!\!
 \frac{(2\pi)^{2\alpha|\setu|}}{\gamma_\setu} \bigg(
 \sum_{\bsm_\setu\le\bsalpha_\setu} (|\bsm_\setu|+1)!\,
 \calS(\bsalpha_\setu,\bsm_\setu)\, \prod_{j\in\setu} b_j^{m_j}\!
 \bigg)^2. \label{eq:final-per}
\end{align}

\item [\textnormal{(c)}] Continuing from \textnormal{(b)}, the norm
    \eqref{eq:norm-K} of both functions satisfy for all components
    $p$,
\begin{align}
 &\max\big\{
 \big\|G(\cdot)_p\big\|_{\calK_{\alpha,\bsgamma}}, 
 \big\|G^{[L]}_\theta(\cdot)_p\big\|_{\calK_{\alpha,\bsgamma}} \hphantom{,}
 \big\} \nonumber\\
 &\le C
 \max_{\setu\subseteq\{1:s\}} \!\!
 \frac{(2\pi)^{\alpha|\setu|}}{\gamma_\setu} \!\!
 \sum_{\bsm_\setu\le\bsalpha_\setu} |\bsm_\setu|!\,
 \calS(\bsalpha_\setu,\bsm_\setu)\, \prod_{j\in\setu} b_j^{m_j},
 \label{eq:G-K}
 \\
 &\big\|\big(G(\cdot)_p - G^{[L]}_\theta(\cdot)_p\big)^2\big\|_{\calK_{\alpha,\bsgamma}} \nonumber \\
 &\le 4\,C^2
 \max_{\setu\subseteq\{1:s\}} \!\!
 \frac{(2\pi)^{\alpha|\setu|}}{\gamma_\setu} \!\!
 \sum_{\bsm_\setu\le\bsalpha_\setu} (|\bsm_\setu|+1)!\,
 \calS(\bsalpha_\setu,\bsm_\setu)\, \prod_{j\in\setu} b_j^{m_j}.
 \label{eq:final-K}
\end{align}
\end{itemize}
\end{theorem}

\begin{proof} With the restrictions
\eqref{eq:demand1}--\eqref{eq:demand3}, it is clear that the right-hand
sides of \eqref{eq:reg-np} and \eqref{eq:reg-per} are bounded from above
by the right-hand sides of \eqref{eq:tar-np} and \eqref{eq:tar-per},
respectively. The bounds \eqref{eq:G-np}, \eqref{eq:G-per} and
\eqref{eq:G-K} then follow from the definition of the norms
\eqref{eq:norm-np}, \eqref{eq:norm-per} and \eqref{eq:norm-K}. The final
bounds \eqref{eq:final-np}, \eqref{eq:final-per} and \eqref{eq:final-K}
can be obtained from \eqref{eq:tar-np} and \eqref{eq:tar-per} using the
Liebniz product rule, see Section~\ref{sec:norm}.
\end{proof}

We now choose weights $\gamma_\setu$ and construct lattice
generating vectors in all three cases to achieve convergence rates that
are independent of the input dimension $s$.

\begin{theorem}[Tailored lattice training points] \label{thm:err}
Under the same conditions as in \textnormal{Theorem~\ref{thm:diff}} and for the
respective settings \textnormal{(a)}, \textnormal{(b)}, \textnormal{(c)}, assume further
that there is a ``summability exponent'' $p^* \in (0,1)$ such that
$\sum_{j\ge 1} b_j^{p^*} < \infty$.
\begin{itemize}
\item [\textnormal{(a)}] Construct as in Subsection~\ref{sec:qmc}$(a)$ a good generating vector $\bsz$ for a randomly-shifted lattice rule in weighted Sobolev space with weights
\begin{align} \label{eq:weight-np}
 &\gamma_\setu := \bigg(
 (|\setu|+1)!\,
 \prod_{j\in\setu} \bigg(
 \sqrt{\frac{(2\pi)^{2\lambda}}{2\zeta(2\lambda)\,2^\lambda}}\, b_j\bigg)
 \bigg)^{\frac{2}{1+\lambda}},
 \\
 &\lambda :=
 \begin{cases}
 \frac{1}{2-2\delta}, \delta\in (0,\frac{1}{2}),
 & \mbox{if } p^* \in (0, \frac{2}{3}], \\
 \frac{p^*}{2-p^*} & \mbox{if } p^* \in (\frac{2}{3}, 1). \nonumber
 \end{cases}
\end{align}
If a non-periodic DNN \eqref{eq:DNN-np} is trained using these lattice
points, and the training error reaches the threshold ${\tt tol} \in
(0,1)$, then the generalization error is bounded by
\begin{align} \label{eq:rate-np}
  \calE_G 
  \le {\tt tol} + \calO\big(N^{-r/2}\big),
  \qquad
  r := \min\Big(1-\delta,\;\frac{1}{p^*}-\frac{1}{2}\Big),
  \quad \delta\in (0,\tfrac{1}{2}).
\end{align}

\item [\textnormal{(b)}] Construct as in Subsection~\ref{sec:qmc}$(b)$ a good generating vector $\bsz$ for  a lattice rule in weighted Hilbert Korobov space with weights
\begin{align} \label{eq:weight-per}
 &\gamma_\setu :=
 (2\pi)^{2\alpha|\setu|}
 \sum_{\bsm_\setu\le\bsalpha_\setu}\!\!
 \bigg(
 (|\bsm_\setu|+1)!\,
 \prod_{j\in\setu}
 \frac{b_j^{m_j}\calS(\alpha,m_j)}{\sqrt{2\zeta(2\alpha\lambda)}}
 \bigg)^{\frac{2}{1+\lambda}},
 \\
 &\alpha := \left\lfloor \tfrac{1}{p^*} + \tfrac{1}{2} \right\rfloor,
 \,
 \lambda := \tfrac{p^*}{2-p^*}. \nonumber
\end{align}
If a periodic DNN \eqref{eq:DNN-per} is trained using these lattice
points, and the training error reaches the threshold ${\tt tol} \in
(0,1)$, then the generalization error is bounded by
\begin{align} \label{eq:rate-per}
  \calE_G 
  \le {\tt tol} + \calO\big(N^{-r/2}\big),
  \qquad
  r := \frac{1}{p^*}-\frac{1}{2}.
\end{align}

\item [\textnormal{(c)}] Construct as in Subsection~\ref{sec:qmc}$(c)$ a good generating vector $\bsz$ for a lattice rule in weighted non-Hilbert Korobov space with weights
\begin{align} \label{eq:weight-K}
 \gamma_\setu :=
 (2\pi)^{\alpha|\setu|}
 \sum_{\bsm_\setu\le\bsalpha_\setu}
 (|\bsm_\setu|+1)!\,
 \prod_{j\in\setu} \big( b_j^{m_j}\calS(\alpha,m_j)\big), \quad
 \alpha := \left\lfloor \tfrac{1}{p^*} \right\rfloor + 1.
\end{align}
If a periodic DNN \eqref{eq:DNN-per} is trained using these lattice
points, and the training error reaches the threshold ${\tt tol} \in
(0,1)$, then the generalization error is bounded by
\begin{align} \label{eq:rate-K}
  \calE_G 
  \le {\tt tol} + \calO\big(N^{-r/2}\big),
  \qquad
  r := \frac{1}{p^*}.
\end{align}
\end{itemize}
In all three settings, the implied constant in the big-$\calO$ bound is
independent of $s$, but depends on the summability exponent $p^*$ via the parameter~$\lambda$ and the weights~$\gamma_\setu$, as well as the arbitrarily small parameter~$\delta$ in setting~\textnormal{(a)} and the smoothness parameter~$\alpha$ in settings~\textnormal{(b)} and~\textnormal{(c)}. 

In practice, if the condition \eqref{eq:demand2} is not achieved after
training, i.e., if $\beta_j \le b_j/S_L$ is not true for some $j$, with $S_L$ defined in \eqref{eq:CSP},
but instead there is a constant $\kappa\ge 1/s_L$ such that
\begin{align} \label{eq:kappa}
  \beta_j \le \kappa\,b_j
  \quad\mbox{for all } j = 1,\ldots,s,
\end{align}
then the bounds \eqref{eq:final-np}, \eqref{eq:final-per},
\eqref{eq:final-K} hold with each $b_j$ replaced by $\kappa\,S_L\,b_j$.
Even though the lattice training points were constructed with the weights
\eqref{eq:weight-np}, \eqref{eq:weight-per}, \eqref{eq:weight-K} based on
the original $b_j$, the big-$\calO$ bounds \eqref{eq:rate-np},
\eqref{eq:rate-per}, \eqref{eq:rate-K} hold now with enlarged implied
constants that depend on $\kappa\,S_L$. In particular, the enlarged constants for \eqref{eq:rate-np} and \eqref{eq:rate-per} are still bounded independently of $s$, but the enlarged constant for \eqref{eq:rate-K} grows exponentially with $s$.
\end{theorem}

\begin{proof}
The results follow from \eqref{eq:combine} by combining
\eqref{eq:final-np} with \eqref{eq:wce-np}, \eqref{eq:final-per} with
\eqref{eq:wce-per}, and \eqref{eq:final-K} with \eqref{eq:wce-K}, and then
choosing the weights and other parameters such as $\alpha$ and $\lambda$
with respect to $p^*$ to ensure that the implied constant is finite in
each case. 
The derivations of weights and convergence rates follow the same arguments as in several recent papers, 
see e.g., \cite{KSS12,KN16} for~(a), \cite{KKKNS22} for~(b), and
\cite{KKS20} for~(c).

The final statement regarding dimension independence for the enlarged constants in the settings (a) and (b) is true because the summability exponent is unchanged, i.e.,
$\sum_{j\ge 1} (\kappa\,S_L\,b_j)^{p^*} = (\kappa\,S_L)^{p^*} \sum_{j\ge
1} b_j^{p^*} < \infty$. The setting (c) is different and will incur an extra factor $(\kappa\,S_L)^{\alpha\,s}$. Detailed derivations for all three settings are given in our survey article \cite{KKNS25b}.
\end{proof}

The theoretical convergence rates above are very conservative because we
were aiming for implied constants that do not depend on $s$. Moreover, the
estimate $|\calE_G - \calE_T| \le \sqrt{|\calE_G^2 - \calE_T^2|}$ also
means we lose half of the rate of convergence, a loss that may not actually be
observed in practice, noting that the first error bound \eqref{eq:err1} also holds.

We see that the convergence rate in setting~(c) is better than in setting~(b), demonstrating an advantage of the non-Hilbert space setting. However, setting~(c) is less forgiving if the condition \eqref{eq:demand2} is not achieved after training, since it will result in a constant that grows exponentially with $s$ while the constant for setting~(b) remains bounded independently of~$s$.

\section{Examples} \label{sec:pde}

\subsection{A parametric diffusion problem}

A specific example is the parametric diffusion problem
\begin{equation}\label{eq:diffusion}
\begin{cases}
- \nabla\cdot \big( a(\bsx, \bsy) \nabla u(\bsx, \bsy)\big)
= f(\bsx), &\bsx \in D, \quad \bsy \in Y = [0,1]^s,\\
u(\bsx,\bsy) = 0, &\bsx \in \partial D,
\end{cases}
\end{equation}
where $D$ is a bounded spatial domain in $\bbR^d$ with $d\in \{1,2,3\}$,
$f\in L^2(D)$, and the diffusion coefficient $a(\bsx,\bsy)$ is a given
positive function of $\bsx$ and $\bsy$. We consider two models:
\begin{itemize}
\item [\textnormal{(a)}] Affine (non-periodic) model (see e.g.,
    \cite{CD15})
\begin{align} \label{eq:affine}
  a(\bsx,\bsy) := \psi_0(\bsx) + \sum_{j=1}^s (y_j - \tfrac{1}{2})\, \psi_j(\bsx).
\end{align}
\item [\textnormal{(b)}] Periodic model (see \cite{KKKNS22,KKS20} but
    without the $1/\sqrt{6}$ scaling)
\begin{align} \label{eq:aper}
  a(\bsx,\bsy) := \psi_0(\bsx) + \sum_{j=1}^s \sin(2\pi y_j)\, \psi_j(\bsx).
\end{align}
\end{itemize}
In both cases we assume that there exist constants $a_{\min}$ and
$a_{\max}$ (different for the affine and periodic models) such that $0 <
a_{\min} \le a(\bsx,\bsy) \le a_{\max} < \infty$, and that the functions
$\psi_j\in L^\infty(D)$ satisfy additional assumptions. An important
assumption when the dimensionality $s$ is high (or infinite) is that there
exists a ``summability exponent'' $p^*\in (0,1)$ such that $\sum_{j\ge 1}
\|\psi_j\|_{\infty}^{p^*} < \infty$.

As a concrete example, we consider a simple observable which is a linear
functional on the PDE solution:
\begin{align} \label{eq:linG}
 G(\bsy) := \calG(u(\cdot,\bsy)) = \int_{D'} u(\bsx, \bsy) \,\rd \bsx,
\end{align}
where $D'$ may be one half/quadrant/octant of the domain $D$. Thus
$N_{\mathrm{obs}}=1$. More generally, the observables might be the vector
of values of the solution $u$ at the points of a finite element mesh, in
which case $N_{\mathrm{obs}}$ could be the number of interior finite
element mesh points.

The space $V = H^1_0(D)$ is the typical Sobolev space with first order
derivatives with respect to $\bsx\in D$, with zero boundary condition and
norm $\|v\|_V := \| \nabla v\|_{L^2(D)}$. The space $V' = H^{-1}(D)$ is
its dual. The parametric regularity of the PDE solutions are known in both cases:
\begin{itemize}
\item [\textnormal{(a)}] For the affine (non-periodic) model
    \eqref{eq:affine} we have (see e.g., \cite{KSS12})
\begin{align} \label{eq:pde-np}
 \|\partial^\bsnu u(\cdot,\bsy)\|_V
 \le \frac{\|f\|_{V'}}{a_{\min}}\, |\bsnu|!\,\bsb^\bsnu, \qquad
 b_j := \frac{\|\psi_j\|_{L_\infty(D)}}{a_{\min}}.
\end{align}
The target function \eqref{eq:linG} satisfies the regularity bound
\eqref{eq:tar-np} with $b_j$ as in \eqref{eq:pde-np} and
\begin{align} \label{eq:C-np}
  C := \frac{\|\calG\|_{V'}\,\|f\|_{V'}}{a_{\min}}.
\end{align}

\item [\textnormal{(b)}] For the periodic model \eqref{eq:aper} we
    have (see \cite{KKS20})
\begin{align} \label{eq:pde-per}
 \|\partial^\bsnu u(\cdot,\bsy)\|_V
 \le \frac{\|f\|_{V'}}{a_{\min}} (2\pi)^{|\bsnu|}
 \sum_{\bsm\le\bsnu} |\bsm|!\,\bsb^\bsm\, \calS(\bsnu,\bsm),
\end{align}
where $b_j$ is defined as in \eqref{eq:pde-np}, but with a different
$a_{\min}$ appropriate for the periodic model. The target function
\eqref{eq:linG} satisfies the regularity bound \eqref{eq:tar-per} with
$b_j$ as in \eqref{eq:pde-np} and $C$ as in \eqref{eq:C-np}, but again
with a different $a_{\min}$. Note that we do not have a $1/\sqrt{6}$
factor in \eqref{eq:aper} unlike in \cite{KKS20,KKKNS22}.
\end{itemize}

\subsection{An algebraic equation}

We consider an algebraic equation $a(\bsy)\,G(\bsy) = 1$, $\bsy\in U_s$,
which mimics the features of the parametric PDE \eqref{eq:diffusion} but
avoids the spatial variable $\bsx$ and the added complexity of a finite
element solver. The target function we seek to approximate by a DNN is
\begin{align} \label{eq:alg-G}
  G(\bsy) := \frac{1}{a(\bsy)},
\end{align}
where we consider two models for the denominator (compare to
\eqref{eq:affine} and \eqref{eq:aper}):
\begin{itemize}
\item [\textnormal{(a)}] Affine (non-periodic) algebraic model
\begin{align} \label{eq:alg-affine}
  a(\bsy) := 1 + \sum_{j=1}^s (y_j - \tfrac{1}{2})\,\psi_j,
  \quad \psi_j := \frac{\eta}{j^q} \mbox{ for } q>1,
  \quad a_{\min} := 1 - \frac{\eta\,\zeta(q)}{2},
\end{align}
where $\zeta(\cdot)$ is again the Riemann zeta function. The target
function \eqref{eq:alg-G} satisfies the regularity bound
\eqref{eq:tar-np} with $b_j = \psi_j/a_{\min}$ and $C = 1/a_{\min}$.

\item [\textnormal{(b)}] Periodic algebraic model
\begin{align} \label{eq:alg-aper}
  a(\bsy) := 1 + \sum_{j=1}^s \sin(2\pi y_j)\,\psi_j,
  \quad \psi_j := \frac{\eta}{j^q} \mbox{ for } q>1,
  \quad a_{\min} := 1 - \eta\,\zeta(q).
\end{align}
The target function \eqref{eq:alg-G} satisfies the regularity bound
\eqref{eq:tar-per} with $b_j = \psi_j/a_{\min}$ and $C = 1/a_{\min}$,
remembering that $a_{\min}$ is different between the affine and
periodic cases.
\end{itemize}
In both cases, the summability exponent $p^*$ in $\sum_{j\ge 1} b_j^{p^*}<\infty$ satisfies $p^*\approx 1/q$. 

\section{Numerical experiments} \label{sec:num}

In our experiments we will train DNNs to approximate the target function
\eqref{eq:alg-G} from the periodic algebraic model \eqref{eq:alg-aper}. We
aim to demonstrate that restricting network parameters to match the
regularity features of the target function as in Theorem~\ref{thm:diff} can
improve the resulting DNN.

\subsection{Tailored regularization} \label{sec:regularization}

We train our DNN to obtain
\begin{align*}
  \theta^* := \argmin_{\theta\in\Theta}
  \bigg(
  \frac{1}{N} \sum_{k=1}^N \big\|
  G(\bsy_k) - G_{\theta}^{[L]}(\bsy_k) \big\|_2^2
  + \lambda\, \|\theta\|_2^2
  + \lambda_1\,\calR_1(\theta)
  \bigg).
\end{align*}
The case $\lambda_1 = 0$ corresponds to the
standard $\ell_2$ regularization. Otherwise we have our tailored
regularization designed based on Theorem~\ref{thm:diff}:
\begin{align}
  \calR_1(\theta) &:=
  \frac{1}{s} \sum_{j=1}^s \frac{1}{d_1} \sum_{p=1}^{d_1} \Big(W_{0,p,j}^2\,
  \frac{L^2}{b_j^2}\Big)^{m/2},
  \label{eq:R1-def} 
\end{align}
with a positive even integer $m = 6$.

The regularization term $\calR_1(\theta)$ is chosen to encourage the
condition \eqref{eq:demand2} from Theorem~\ref{thm:diff}, i.e.,
$\|W_{0,:,j}\|_\infty \le\beta_j \le b_j/S_L$ for all $j=1,\ldots,s$,
which is equivalent to
\[
  \max_{1\le j\le s} \max_{1\le p\le d_1} |W_{0,p,j}|\, \frac{S_L}{b_j} \le 1.
\]
In an ideal scenario we could have $\xi = \tau = 1$ in
\eqref{eq:common} and $\rho = 1$ in \eqref{eq:demand1}, leading to $S_L\le
L$ in~\eqref{eq:demand1b}. While this is unlikely to be true in general,
for simplicity we replace $S_L$ by $L$ in designing our regularization
term, to at least capture the fact that $S_L$ will generally increase with
increasing~$L$. Ultimately this only affects the scaling of the
regularization term.

Recall that $\max_{1\le i\le n} |a_i| = \lim_{m\to\infty}
(\frac{1}{n}\sum_{i=1}^n |a_i|^m)^{1/m}$. Since we need $\calR_1(\theta)$
to be smooth so that an automatic differentiation routine can compute
$\nabla_\theta\calR_1(\theta)$, we replace the two maximums by the average
of the terms raised to an even power $m = 6$ as follows:
\[
  \bigg[\frac{1}{s} \sum_{j=1}^s \bigg(\bigg[
  \frac{1}{d_1} \sum_{p=1}^{d_1} \bigg(|W_{0,p,j}|\, \frac{L}{b_j}\bigg)^m\bigg]^{1/m}
  \bigg)^m \bigg]^{1/m} \le 1,
\]
which simplifies to
\[
  \frac{1}{s} \sum_{j=1}^s
  \frac{1}{d_1} \sum_{p=1}^{d_1} \bigg(W_{0,p,j}^2\, \frac{L^2}{b_j^2}\bigg)^{m/2}
  \le 1^m = 1.
\]
The even power allows us to omit the absolute values. We note that the
function $x^6$ stays mostly constant and below $1$ for $x\le 1$ and grows
very quickly for $x>1$.

The conditions \eqref{eq:demand1} and \eqref{eq:demand3} restrict the magnitude of $R_1,\ldots,R_{L-1}$ and~$R_L$. These are already controlled by the standard $\ell_2$ regularization term $\lambda\,\|\theta\|_2^2$.
The condition \eqref{eq:demand2} is the most important of the three conditions, since the rate of decay of the sequence $\beta_j$ will affect the empirical convergence rate for QMC methods and has a direct impact on whether or not the implied constant grows exponentially with~$s$. 

After training, from our trained DNN we can compute $\beta_1,\ldots,\beta_s$, and thus obtain a value of $\kappa$ satisfying \eqref{eq:kappa}. We can also compute $R_1,\ldots,R_{L-1}$, thus obtain a value of $\rho$ satisfying \eqref{eq:demand1}. Furthermore, we can compute $R_L$, $C_L$ and $S_L$ defined in \eqref{eq:CSP}. We see from \eqref{eq:demand1b} that $S_L$ will generally scale like~$\rho^L$.

In practice, if $C_L\le C$ is \emph{not} true as demanded in \eqref{eq:demand3}, then we can adjust our norm bounds in Theorem~\ref{thm:diff} with $C$ replaced by $C_L$. If $\kappa \le 1/S_L$ is \emph{not} true, i.e., $\beta_j \le b_j/S_L$ for all $j$ is \emph{not} true as demanded in \eqref{eq:demand2} (and most likely it will not be true), then we can adjust our norm bounds in Theorem~\ref{thm:diff} with each $b_j$ replaced by $\kappa\,S_L\,b_j$. Consequently, as indicated at the end of Theorem~\ref{thm:err}, the generalization error bounds still hold, but with enlarged implied constants dependent on the value of~$\kappa\,S_L$.

\subsection{Training with full-batch Adam algorithm}

We utilize the Python-based open-source deep learning framework PyTorch \cite{PyTorch}
to execute a full-batch Adam (adaptive moment estimation) optimization algorithm \cite{Adam}
which is a gradient descent method. The algorithm begins with
a specified initial parameter configuration (random ``Glorot initialization'') and employs
a ``learning rate'' (step size) of $10^{-4}$. The optimization process halts either
when the training error $\calE_T$ (see \eqref{eq:ET}) reaches the
threshold ${\tt tol} = 10^{-3}$, or upon completing a maximum number of ``epochs'' (steps)
of $40000$. 

We use a randomly shifted version of a lattice point set as our training
points $\{\bsy_k\}_{k\ge 1}$, and progressively train the network using $N = 2^5, \ldots,
2^{12}$ points. For simplicity, we take one ``off-the-shelf'' embedded lattice rule generating vector, 
constructed according to the method described in \cite{CKN06}. 

We illustrate our results for the case $N_{\rm obs} = 1$. We compare the training error $\calE_T$ to an estimate of the generalization error $\calE_G$ (see \eqref{eq:EG}) given
by
\[
  \calE_G
  \approx
  \widetilde{\calE}_G
  := \bigg(\frac{1}{M} \sum_{i=1}^M \big(G(\bst_i) - G_{\theta}^{[L]}(\bst_i) \big)^2 \bigg)^{1/2},
\]
where we use one different and independent random shift with the same embedded lattice rule with a much higher number of points $M = 2^{15} \gg N$. We then compute an estimate of the generalization gap by
$|\widetilde{\calE}_G - \calE_T|$. 
Recall that $|\calE_G^2 - \calE_T^2|$ is the quadrature error for the
integrand $(G - G_\theta^{L})^2$, so we expect to observe a bound
of the form
\begin{align*}
  \calE_G
  \le \calE_T + |\calE_G - \calE_T|
  &\le \calE_T + \sqrt{e^{\rm wor}_N(\bsz,\calW_{\alpha,\bsgamma})\, \| (G - G_{\theta}^{[L]})^2 \|_{\calW_{\alpha,\bsgamma}}}
  \\
  &\le 10^{-3} + \calO\big(N^{-r/2}\big),
\end{align*}
with $r$ depending on the three settings in Theorem~\ref{thm:err}.
Ideally, we hope to observe that the estimated generalization error $\widetilde{\calE}_G$ is consistently reduced when tailored regularization is applied.
Additionally, we hope to see a higher convergence rate for the generalization gap $|\widetilde{\calE}_G - \calE_T|$ when tailored regularization is employed and QMC points are used instead of
standard Monte Carlo (MC) points.

For the numerical experiments, we will investigate the periodic algebraic model
\eqref{eq:alg-G} and \eqref{eq:alg-aper} with
\[
  \eta = 0.5, \quad  q = 2.5.
\]
Recall that the summability exponent satisfies $1/p^* \approx q = 2.5$. The theoretical convergence rate for the generalization error according to Theorem~\ref{thm:err} is $\calO(N^{-r/2})$, with $r = 1/p^*-1/2 \approx 2$ in setting~(b) or with $r = 1/p^* \approx 2.5$ in setting~(c), yielding around the rate $\calO(N^{-1})$ or $\calO(N^{-1.25})$, respectively.

We will train a periodic DNN \eqref{eq:DNN-per} with the sigmoid activation function using two sets of hyperparameters:
\begin{enumerate}
\item 
  $L = 3$,\;\; $N_{\rm obs} = 1 = d_{L+1} = d_4$, \;\; $d_3 = d_2 = d_1 = 32$, $d_0 = s = 50$\; ($3777$~parameters);
\item 
  $L = 12$, $N_{\rm obs} = 1 = d_{L+1} = d_{13}$, $d_{12} = \cdots = d_1 = 30$, $d_0 = s = 50$\;  ($11791$~parameters);
\end{enumerate}
that is, one relatively small network and one larger network.
While some may argue that these networks are small, we point out that ``tiny'' neural networks have recently been shown to be very successful in computer graphics, see \cite{MESK22,MRNK21,VSWAEL23}.

\begin{figure}[t]
\begin{center}
 \includegraphics[width=12.9cm]{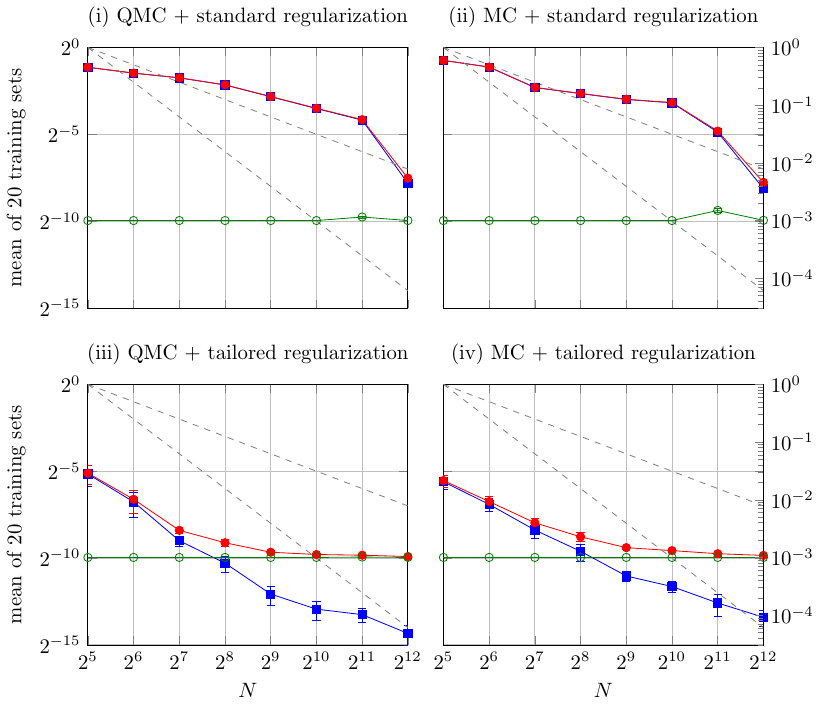}
 \caption{Values of $\calE_T$ (\textcolor{green!50!black}{$\circ$} green circles), $\widetilde{\calE}_G$ (\textcolor{red}{$\bullet$} red dots), and $|\widetilde{\calE}_G - \calE_T|$ (\textcolor{blue}{\scalebox{0.7}{$\blacksquare$}} blue squares) as $N$ increases, with 
 $L = 3$, $N_{\rm obs} = 1$, $d_\ell = 32$, $s = 50$.}
 \label{fig1a}
\end{center}
\end{figure}

\begin{figure}[t]
\begin{center}
 \includegraphics[width=12.9cm]{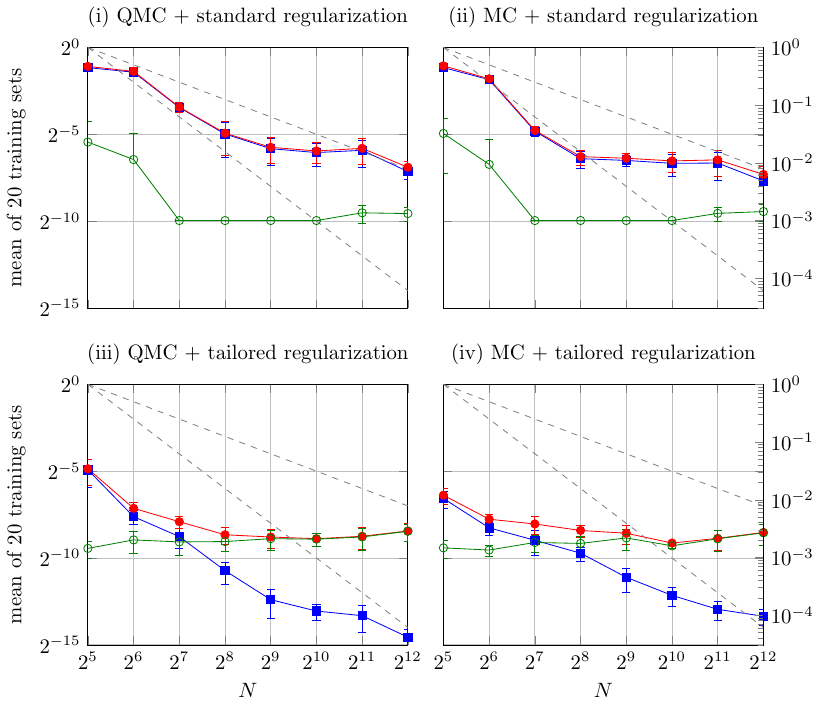}
 \caption{Values of $\calE_T$ (\textcolor{green!50!black}{$\circ$} green circles), $\widetilde{\calE}_G$ (\textcolor{red}{$\bullet$} red dots), and $|\widetilde{\calE}_G - \calE_T|$ (\textcolor{blue}{\scalebox{0.7}{$\blacksquare$}} blue squares) as $N$ increases, with  
 $L = 12$, $N_{\rm obs} = 1$, $d_\ell = 30$, $s = 50$.}
 \label{fig1b}
\end{center}
\end{figure}

In Figure~\ref{fig1a}, for the first set of hyperparameters ($L=3$), we plot 
the values of the training error $\calE_T$~(\textcolor{green!50!black}{$\circ$} green circles), the estimated generalization error $\widetilde{\calE}_G$~(\textcolor{red}{$\bullet$} red dots), and the generalization gap $|\widetilde{\calE}_G - \calE_T|$~(\textcolor{blue}{\scalebox{0.7}{$\blacksquare$}} blue squares) for $N =
2^5,\ldots,2^{12}$, where each of the three quantities is averaged over $20$
different random ``Glorot'' initializations of the network and random shifts of the lattice rule. 
The top two graphs (labelled ``standard regularization'') were obtained with regularization parameter $\lambda = 10^{-8}$ but without our tailored
regularization (so $\lambda_1= 0$), while the bottom two graphs (labelled ``tailored regularization'') were obtained with $\lambda = 10^{-8}$ and with our tailored regularization $\lambda_1 = 10^{-8}$. 
The two graphs on the left (labelled ``QMC'') were obtained by training with
randomly shifted embedded lattice points as mentioned above, while the two graphs on the right (labelled ``MC'') were obtained with uniformly distributed random points for comparison. 
In Figure~\ref{fig1b} we give the same plots for the second set of hyperparameters ($L=12$).

We observe in all four cases in Figure~\ref{fig1a} ($L=3$) that the training error $\calE_T$~(\textcolor{green!50!black}{$\circ$}) is able to reach the threshold ${\tt tol}=10^{-3}$ for most values of $N$. 
With just the standard regularization (i)--(ii), both $\widetilde{\calE}_G$~(\textcolor{red}{$\bullet$}) and $|\widetilde{\calE}_G - \calE_T|$~(\textcolor{blue}{\scalebox{0.7}{$\blacksquare$}}) are large and essentially overlapping, and plateau initially before diving down from $N=2^{11}$, but even at $N=2^{12}$ they are still some distance away from ${\tt tol}=10^{-3}$.
On the other hand, with tailored regularization (iii)--(iv), $|\widetilde{\calE}_G - \calE_T|$~(\textcolor{blue}{\scalebox{0.7}{$\blacksquare$}}) decays at a rate between $\calO(N^{-1})$ and $\calO(N^{-2})$, and quickly crosses below ${\tt tol} =10^{-3}$ by $N=2^9$, with $\widetilde{\calE}_G$~(\textcolor{red}{$\bullet$}) fairly close to ${\tt tol} = 10^{-3}$ from $N = 2^{10}$. 

In Figure~\ref{fig1b} ($L=12$) we observe that, with just the standard regularization (i)--(ii), the training error $\calE_T$~(\textcolor{green!50!black}{$\circ$}) cannot reach the threshold ${\tt tol}=10^{-3}$ until $N=2^7$, and both $\widetilde{\calE}_G$~(\textcolor{red}{$\bullet$}) and $|\widetilde{\calE}_G - \calE_T|$~(\textcolor{blue}{\scalebox{0.7}{$\blacksquare$}}) stay some distance away from ${\tt tol}=10^{-3}$ even at $N=2^{12}$. 
On the other hand, with tailored regularization (iii)--(iv), the training error $\calE_T$~(\textcolor{green!50!black}{$\circ$}) cannot reach the threshold ${\tt tol}=10^{-3}$ for all values of $N$. While this may be disappointing, we observe that $\widetilde{\calE}_G$~(\textcolor{red}{$\bullet$}) follows $\calE_T$~(\textcolor{green!50!black}{$\circ$}) closely, and the gap $|\widetilde{\calE}_G - \calE_T|$~(\textcolor{blue}{\scalebox{0.7}{$\blacksquare$}}) crosses below $\calE_T$~(\textcolor{green!50!black}{$\circ$}) from $N = 2^7$ and decays at a rate between $\calO(N^{-1})$ and $\calO(N^{-2})$ as before. Note that the learning rate is constant. We expect the results to improve once we introduce a strategy to reduce the learning rate during training.

The empirical convergence rate between $\calO(N^{-1})$ and $\calO(N^{-2})$ is consistent with our theoretical prediction of around $\calO(N^{-1})$ or $\calO(N^{-1.25})$ for the case $1/p^*\approx q = 2.5$ in Theorem~\ref{thm:err}, settings (b) or (c). However, since we used an off-the-shelf lattice generating vector rather than one explicitly constructed for either setting, this prediction should be regarded as an indication rather than a guarantee. We also carried out experiments with a range of values of $q$. Smaller values of $q$ gave lower convergence rates but larger values of $q$ did not always conclusively improve the rates. We can speculate about the reasons, e.g., the theoretical implied constant is larger in the higher smoothness setting and we need larger values of $N$ in order to observe the improved rate in practice. Furthermore, we did not attempt to optimize DNN training and instead employed a simple and easily reproducible training procedure. Experiments with explicitly constructed lattice rules for different values of~$q$ and more advanced training methods are left for future research.

We also experimented with a non-periodic DNN with lattice training points on the periodic algebraic function. Not surprisingly, we observed lower convergence rate; this is consistent with Theorem~\ref{thm:err}(a) which predicts close to $\calO(N^{-1/2})$ convergence in this setting. To obtain higher order convergence rate with a non-periodic DNN, we would need to tap into different QMC theory, e.g., using tent-transformed lattice rules or higher order digital nets, building on the parametric regularity bounds in this paper.

In summary, tailored regularization provides consistent improvement in all our experiments, and it benefits QMC as well as MC. A surprise is that MC performs just as well as QMC, with QMC perhaps having a slightly better empirical convergence rate for the gap $|\widetilde{\calE}_G - \calE_T|$~(\textcolor{blue}{\scalebox{0.7}{$\blacksquare$}}). The theory in our paper does not directly apply to MC, although loosely speaking we may replace the worst case error bound $e^{\rm wor}_N(\bsz, \calW_{\alpha,\bsgamma})\, \| (G - G_{\theta}^{L})^2\|_{\calW_{\alpha,\bsgamma}}$ by the variance of $(G - G_{\theta}^{L})^2$ divided by~$\sqrt{N}$. A possible explanation for the success of MC with tailored regularization may be that tailored regularization reduces the variance of the squared difference.

\begin{figure}[t]
\begin{center}
 \includegraphics[width=11cm]{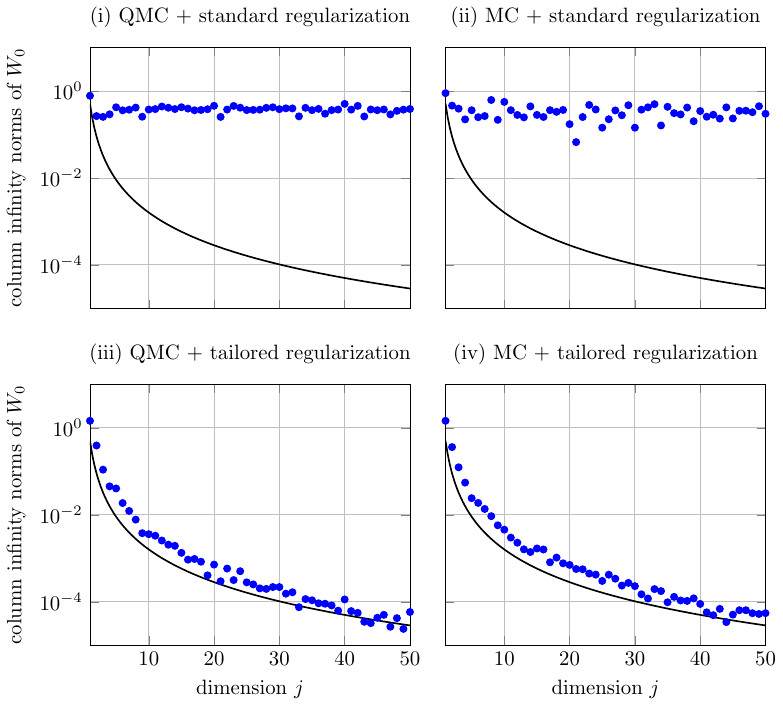}
 \caption{Values of $\log(\beta_j)$ (\textcolor{blue}{$\bullet$} blue dots) 
 and $\log(b_j/L)$ (black line) for $j=1,\ldots,s$ for one random initialization and one random shift, with 
 $L = 3$, $N_{\rm obs} = 1$, $d_\ell = 32$, $s = 50$, $N=2^5$.}
 \label{fig2a}
\end{center}
\end{figure}

\begin{figure}[t]
\begin{center} 
 \includegraphics[width=11cm]{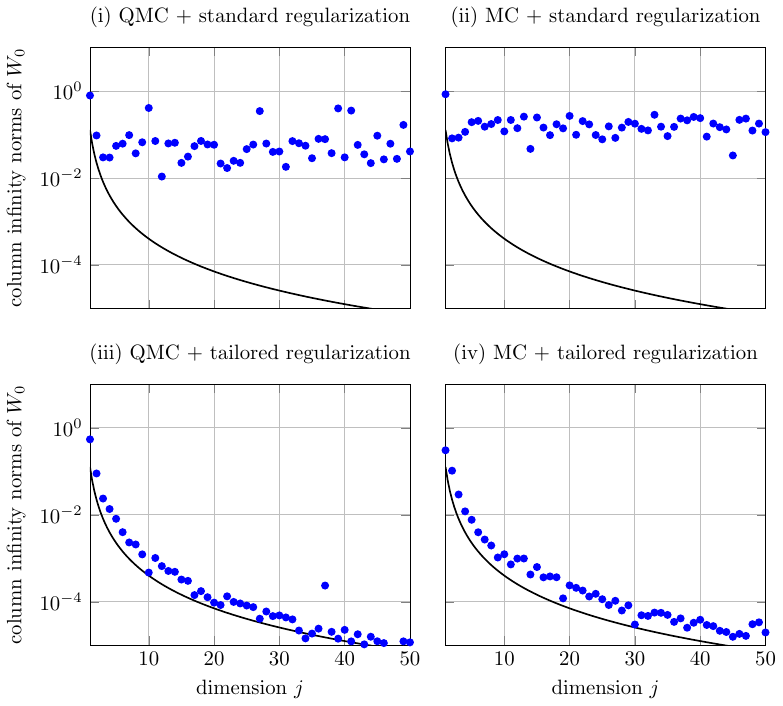}
 \caption{Values of $\log(\beta_j)$ (\textcolor{blue}{$\bullet$} blue dots) 
 and $\log(b_j/L)$ (black line) for $j=1,\ldots,s$, for one random initialization and one random shift, with 
 $L = 12$, $N_{\rm obs}\! =\! 1$, $d_\ell\! =\! 30$, $s = 50$, $N=2^5$.}
 \label{fig2b}
\end{center}
\end{figure}

In Figure~\ref{fig2a}, for one random initialization and one random shift with the first set of hyperparameters ($L=3$) and $N=2^5$, we plot the values of $\log(\beta_j)$ as blue dots and $\log(b_j/L)$ as a black line
for $j=1,\ldots,s$. We see that the sequence $(\beta_j)$ with tailored regularization (iii)--(iv) decays closely following the sequence $(b_j/L)$ but there is a scaling
constant. In fact, for this instance of case~(iii) we achieve $\beta_j\le \kappa\,b_j$ for all $j=1,\ldots,s$ with $\kappa = 1.4851$ and for the matrices $W_1,\ldots, W_{L-1}$ we achieve $R_\ell\le \rho$ for all $\ell=1,\ldots,L-1$ with $\rho = 18.371$. Results for the second set of hyperparameters ($L=12$) are shown in Figure~\ref{fig2b}, exhibiting a similar trend, where for this instance of case~(iii) we achieve $\kappa = 1.2951$ and $\rho = 10.997$.

\section{Technical proofs} \label{sec:proof}

\subsection{Proof of regularity bound for non-periodic DNN} \label{sec:reg-np}

Recall the notation from Section~\ref{sec:reg}, including $\bbN_0 :=
\{0,1,2,\ldots\}$ and $\indx := \{\bsnu\in\bbN_0^\infty :
|\bsnu|<\infty\}$. In the following we will use the Kronecker delta
function: $\delta_{\bsnu,\bszero}$ is $1$ if $\bsnu=\bszero$ and is $0$
otherwise. Let $\bse_j$ denote the multiindex whose $j$th component is
$1$ and all other components are $0$.

First we prove a generic technical result involving a sharp recursion. The
key point is that this result separates the ``product'' (see $\beta_j$
below) contributions from other ``order dependent'' (see $\Gamma_n$
below) contributions of the sequence. 

\begin{lemma} \label{lem:recur-np}
Let the nonnegative sequence
$(\bbA_{\bsnu,\lambda})_{\bsnu\in\indx,\,\lambda\in\bbN_0}$ satisfy
$\bbA_{\bsnu,0} = \delta_{\bsnu,\bszero}$, $\bbA_{\bsnu,\lambda} = 0$ for
$\lambda > |\bsnu|$, and otherwise satisfy the recursion
\begin{align} \label{eq:A1-np}
  \bbA_{\bsnu+\bse_j,\lambda}
  = R \sum_{\bsm\le\bsnu} \binom{\bsnu}{\bsm}
  \bsbeta^{\bsnu-\bsm+\bse_j}\, \Gamma_{|\bsnu-\bsm|+1}\, \bbA_{\bsm,\lambda-1},
  \quad \bsnu\in\indx,
  \; j\ge 1,
  \; 1\le \lambda\le |\bsnu|+1,
\end{align}
for some positive $R$, $(\beta_j)_{j\ge 1}$ and $(\Gamma_n)_{n\ge 0}$. Then for all
$\bsnu\ne\bszero$ and $1\le \lambda\le |\bsnu|$ we have
\begin{align} \label{eq:A2-np}
  \bbA_{\bsnu,\lambda}
  = R^\lambda\, \bsbeta^\bsnu\, \bbB_{|\bsnu|,\lambda}\,,
\end{align}
with
\begin{align} \label{eq:Bn}
 \bbB_{n,1} := \Gamma_n \quad\mbox{for $n\ge 1$}, \quad
 \bbB_{n,\lambda} :=
 \sum_{i=\lambda-1}^{n-1} \binom{n-1}{i} \Gamma_{n-i}\,\bbB_{i,\lambda-1}
 \quad\mbox{for $n\ge\lambda\ge 2$}.
\end{align}
If \eqref{eq:A1-np} is replaced by an inequality $\le$ then so is
\eqref{eq:A2-np}.
\end{lemma}

\begin{proof}
We prove \eqref{eq:A2-np} by induction on $\lambda$. Consider first the case
$\lambda = 1$. Then for all $\bsnu\in\indx$ and all $j\ge 1$
we have from \eqref{eq:A1-np} (since only the
$\bsm = \bszero$ term remains)
\[
  \bbA_{\bsnu+\bse_j,1}
  = R\, \bsbeta^{\bsnu+\bse_j}\,\Gamma_{|\bsnu|+1}
  = R\, \bsbeta^{\bsnu+\bse_j}\,\bbB_{|\bsnu|+1,1},
\]
where we used $\Gamma_{|\bsnu|+1} = \bbB_{|\bsnu|+1,1}$ from
\eqref{eq:Bn}. This proves that \eqref{eq:A2-np} holds for $\lambda =
1$.

For $\lambda \ge 2$, we assume as induction hypothesis that
(c.f.~\eqref{eq:A2-np})
\[
  \bbA_{\bsm,\lambda-1} = R^{\lambda-1}\,
  \bsbeta^\bsm\,\bbB_{|\bsm|,\lambda-1}
\]
for all $\bsm\ne\bszero$ with $|\bsm|\ge \lambda-1$. Recall that
$\bbA_{\bsm,\lambda-1}=0$ for $|\bsm|< \lambda-1$. Then for all $\bsnu\in\indx$ and all $j\ge 1$
we have from \eqref{eq:A1-np} that
\begin{align*}
  \bbA_{\bsnu+\bse_j,\lambda}
  &= R
  \sum_{\substack{\bsm\le\bsnu \\ |\bsm|\ge\lambda-1}} \binom{\bsnu}{\bsm}
  \bsbeta^{\bsnu-\bsm+\bse_j}\,\Gamma_{|\bsnu-\bsm|+1}\,
  R^{\lambda-1}\, 
  \bsbeta^\bsm\,\bbB_{|\bsm|,\lambda-1}  \\
  &= R^\lambda\, 
  \bsbeta^{\bsnu+\bse_j}
  \sum_{\substack{\bsm\le\bsnu \\ |\bsm|\ge\lambda-1}} \binom{\bsnu}{\bsm}
  \Gamma_{|\bsnu-\bsm|+1} \, \bbB_{|\bsm|,\lambda-1}\,  \\
  &= R^\lambda\, 
  \bsbeta^{\bsnu+\bse_j}
  \sum_{i=\lambda-1}^{|\bsnu|}
  \sum_{\substack{\bsm\le\bsnu \\ |\bsm|=i}} \binom{\bsnu}{\bsm}
  \Gamma_{|\bsnu|-i+1} \,
  \bbB_{i,\lambda-1} \\
  &= R^\lambda\, 
  \bsbeta^{\bsnu+\bse_j}
  \underbrace{\sum_{i=\lambda-1}^{|\bsnu|}
  \binom{|\bsnu|}{i} \Gamma_{|\bsnu|-i+1} \,
  \bbB_{i,\lambda-1}}_{=\,\bbB_{|\bsnu|+1,\lambda}},
\end{align*}
where we used the combinatorial identity (e.g., consider the number of
ways to select $i$ distinct balls from some baskets containing a total
number of $|\bsnu|$ distinct balls)
\begin{align} \label{eq:comb}
 \sum_{\substack{\bsm\le\bsnu \\ |\bsm|=i}} \binom{\bsnu}{\bsm}
 = \binom{|\bsnu|}{i},
\end{align}
together with the recursive definition \eqref{eq:Bn}. This completes the
induction proof.

If the equality in \eqref{eq:A1-np} is replaced by an inequality $\le$\,, then it is easy
to see from the proof that this will lead to an inequality $\le$ in the final
result \eqref{eq:A2-np}.
\end{proof}

\begin{proofof}{Theorem~\ref{thm:der}(a)} Consider first a shallow
network with one hidden layer
\[
  G_\theta^{[1]}(\bsy) = W_1\, \sigma (W_0\,\bsy + \bsv_0) + \bsv_1
  = \bigg[
    \sum_{q=1}^{d_1} W_{1,p,q}\, \sigma\bigg( \sum_{r=1}^{d_0} W_{0,q,r}\, y_r
  + v_{0,q}\bigg) + v_{1,p} \bigg]_{p=1}^{d_2},
\]
where $d_0 = s$. For any $\bsnu\ne\bszero$ and any index $p$, it is easy
to see that by the chain rule every derivative with respect to a variable $y_j$ will yield
an extra constant factor $W_{0,q,j}$, leading to
\[
  \partial^\bsnu G_\theta^{[1]}(\bsy)_p
  = \sum_{q=1}^{d_1} W_{1,p,q}\,  \sigma^{(|\bsnu|)} \bigg( \sum_{r=1}^{d_0} W_{0,q,r}\, y_r
  + v_{0,q}\bigg)\,
  \prod_{j=1}^{d_0} W_{0,q,j}^{\nu_j},
\]
and thus
\begin{align*}
  |\partial^\bsnu G_\theta^{[1]}(\bsy)_p|
  &\le \sum_{q=1}^{d_1} |W_{1,p,q}|\, A_{|\bsnu|}
  \prod_{j=1}^{d_0} |W_{0,q,j}|^{\nu_j} \\
  &\le R_1\, A_{|\bsnu|}  \prod_{j=1}^s \beta_j^{\nu_j}
  = R_1\, \bsbeta^\bsnu\, A_{|\bsnu|}
  = R_1\, \bsbeta^\bsnu\, \Gamma_{|\bsnu|}^{[1]}\,,
\end{align*}
where we used the assumptions \eqref{eq:beta}, \eqref{eq:R},
\eqref{eq:sigma}, together with $A_{|\bsnu|} = \Gamma_{|\bsnu|}^{[1]}$
from \eqref{eq:G-def}. Hence \eqref{eq:der-np} holds for $\ell=1$.

Now consider $\ell \ge 2$ and suppose the result \eqref{eq:der-np} holds
for depth $\ell-1$. We can write
\begin{align*}
  G_\theta^{[\ell]}(\bsy) = W_\ell\, \sigma (G_\theta^{[\ell-1]}(\bsy)) + \bsv_\ell
  = \bigg[ \sum_{q=1}^{d_\ell} W_{\ell,p,q}\,
  \sigma (G_\theta^{[\ell-1]}(\bsy)_q) + v_{\ell,p} \bigg]_{p=1}^{d_{\ell+1}}.
\end{align*}
For any $\bsnu\ne\bszero$ and any index $p$, we have
\begin{align} \label{eq:der2a}
  \partial^\bsnu G_\theta^{[\ell]}(\bsy)_p
  &= \sum_{q=1}^{d_\ell} W_{\ell,p,q}\,
 \partial^\bsnu \big( \sigma (G_\theta^{[\ell-1]}(\bsy)_q)\big).
\end{align}
For each index $q$, the recursive form of the multivariate Fa\`a di Bruno
formula gives (see \cite[formula (3.1) and (3.5)]{Savits06} with the
dimensionality $m$ there set to $1$)
\begin{align} \label{eq:der2b}
  \partial^\bsnu \big( \sigma (G_\theta^{[\ell-1]}(\bsy)_q)\big)
  &=
  \sum_{\lambda=1}^{|\bsnu|} \sigma^{(\lambda)} (G_\theta^{[\ell-1]}(\bsy)_q)\,
  \alpha_{\bsnu,\lambda}(\bsy),
\end{align}
where the auxiliary functions $\alpha_{\bsnu,\lambda}(\bsy)$ 
are defined recursively by $\alpha_{\bsnu,0}(\bsy) :=
\delta_{\bsnu,\bszero}$, $\alpha_{\bsnu,\lambda}(\bsy) := 0$
for $\lambda > |\bsnu|$, and otherwise 
\begin{align} \label{eq:tri-np}
  \alpha_{\bsnu+\bse_j,\lambda}(\bsy) :=
  \sum_{\bsm\le\bsnu} \binom{\bsnu}{\bsm}
  \partial^{\bsnu-\bsm+\bse_j} G_\theta^{[\ell-1]}(\bsy)_q \,
  \alpha_{\bsm,\lambda-1}(\bsy).
\end{align}
Using the induction hypothesis \eqref{eq:der-np} for depth $\ell-1$, we obtain  
\begin{align*} 
  |\alpha_{\bsnu+\bse_j,\lambda}(\bsy)|
  &\le \sum_{\bsm\le\bsnu} \binom{\bsnu}{\bsm}
  |\partial^{\bsnu-\bsm+\bse_j} G_\theta^{[\ell-1]}(\bsy)_q| \,
  |\alpha_{\bsm,\lambda-1}(\bsy)| \\
  &\le \sum_{\bsm\le\bsnu } \binom{\bsnu}{\bsm}
  \Big(
  R_{\ell-1}\, \bsbeta^{\bsnu-\bsm+\bse_j}\,\Gamma_{|\bsnu-\bsm|+1}^{[\ell-1]} \Big)\,
  |\alpha_{\bsm,\lambda-1}(\bsy)|.  
\end{align*}
Taking $R = R_{\ell-1}$ and $\Gamma_n = \Gamma_n^{[\ell-1]}$ in Lemma~\ref{lem:recur-np} (with inequalities $\le$) then yields for all $\bsnu\ne\bszero$,
\begin{align} \label{eq:alpha-np}
  |\alpha_{\bsnu,\lambda}(\bsy)|
  \le R_{\ell-1}^\lambda\,\bsbeta^{\bsnu}\,\bbB_{|\bsnu|,\lambda}^{[\ell-1]},
\end{align}
with the sequence $\bbB_{n,\lambda}^{[\ell-1]}$ as defined in \eqref{eq:B-def}.
Applying \eqref{eq:sigma} and \eqref{eq:alpha-np} in \eqref{eq:der2b} yields
\begin{align} \label{eq:der2c}
  \big|\partial^\bsnu \big( \sigma (G_\theta^{[\ell-1]}(\bsy)_q)\big)\big| 
  &\le \sum_{\lambda=1}^{|\bsnu|} A_\lambda\,
  \Big( R_{\ell-1}^\lambda\,\bsbeta^{\bsnu}\,\bbB_{|\bsnu|,\lambda}^{[\ell-1]} \Big)
  = \bsbeta^{\bsnu}\,\Gamma_{|\bsnu|}^{[\ell]},
\end{align}
with the sequence $\Gamma_n^{[\ell]}$ as defined in \eqref{eq:G-def}.
Using \eqref{eq:R} and \eqref{eq:der2c} in \eqref{eq:der2a}, we arrive at
\begin{align*} 
  \big|\partial^\bsnu G_\theta^{[\ell]}(\bsy)_p\big|
  \le \sum_{q=1}^{d_\ell} |W_{\ell,p,q}|\, \bsbeta^{\bsnu}\,\Gamma_{|\bsnu|}^{[\ell]}
  \le R_\ell\,\bsbeta^{\bsnu}\,\Gamma_{|\bsnu|}^{[\ell]}.
\end{align*}
This completes the proof.
\end{proofof}

The bound in Theorem~\ref{thm:der}(a) arises from bounding the matrix elements \eqref{eq:beta} and \eqref{eq:R}, the derivatives of the activation function \eqref{eq:sigma}, and from taking absolute values with the triangle inequality in \eqref{eq:der2a}--\eqref{eq:tri-np}. Otherwise the induction argument from
Lemma~\ref{lem:recur-np} is sharp.

\subsection{Proof of regularity bound for periodic DNN} \label{sec:reg-per}

As in the non-periodic case, we prove here a generic technical result
involving a sharp recursion for the periodic case. The key point is that
this result separates the ``product'' contributions from other
``order dependent'' contributions of the sequence, as well as the
dependence on Stirling numbers of the second kind. 

\begin{lemma} \label{lem:recur-per}
Let the nonnegative sequence
$(\bbA_{\bsnu,\lambda})_{\bsnu\in\indx,\,\lambda\in\bbN_0}$ satisfy
$\bbA_{\bsnu,0} = \delta_{\bsnu,\bszero}$, $\bbA_{\bsnu,\lambda} = 0$ for
$\lambda > |\bsnu|$, and otherwise satisfy the recursion
\begin{align} \label{eq:A1-per}
  \bbA_{\bsnu+\bse_j,\lambda}
  = R
  \sum_{\bsm\le\bsnu} \binom{\bsnu}{\bsm}
  \bigg(\sum_{\bsw\le\bsnu-\bsm+\bse_j} \bsbeta^{\bsw}\,
  \Gamma_{|\bsw|}\,
  & \calS(\bsnu-\bsm+\bse_j,\bsw)\bigg) \bbA_{\bsm,\lambda-1}, \nonumber\\
  &\qquad \bsnu\in\indx,\; j\ge 1, \; 1\le\lambda\le |\bsnu|+1,
\end{align}
for some positive $R$, $(\beta_j)_{j\ge 1}$ and $(\Gamma_n)_{n\ge 0}$. Then for all
$\bsnu\ne\bszero$ and $1\le \lambda\le |\bsnu|$ we have
\begin{align} \label{eq:A2-per}
  \bbA_{\bsnu,\lambda}
  = R^\lambda
  \sum_{\satop{\bsu\le\bsnu}{|\bsu|\ge\lambda}} \bsbeta^\bsu\,\bbB_{|\bsu|,\lambda}\,
  \calS(\bsnu,\bsu),
\end{align}
with $\bbB_{n,\lambda}$ defined as in \eqref{eq:Bn} but with $\Gamma_n$ as used in \eqref{eq:A1-per} rather than in \eqref{eq:A1-np}.
If \eqref{eq:A1-per} is replaced by an inequality $\le$ then so is
\eqref{eq:A2-per}.
\end{lemma}

\begin{proof}
We prove \eqref{eq:A2-per} by induction on $\lambda$. Consider first the case
$\lambda = 1$. Then for all $\bsnu\in\indx$ and all $j\ge 1$ we have from \eqref{eq:A1-per}
(since only the $\bsm = \bszero$ term remains)
\begin{align*}
  \bbA_{\bsnu+\bse_j,1}
  = R
  \sum_{\bsw\le\bsnu+\bse_j} \bsbeta^\bsw\,\Gamma_{|\bsw|}\,\calS(\bsnu+\bse_j,\bsw)
  = R
  \sum_{\satop{\bsw\le\bsnu+\bse_j}{|\bsw|\ge 1}} 
  \bsbeta^\bsw\, \bbB_{|\bsw|,1}\,\calS(\bsnu+\bse_j,\bsw),
\end{align*}
which agrees with \eqref{eq:A2-per} since $\Gamma_{|\bsw|} =
\bbB_{|\bsw|,1}$ as defined in \eqref{eq:Bn} and $\calS(\bsnu+\bse_j,\bszero) = 0$. This proves that
\eqref{eq:A2-per} holds for $\lambda =1$.

For $\lambda\ge 2$, we assume as induction hypothesis that
(c.f.~\eqref{eq:A2-per})
\[
  \bbA_{\bsm,\lambda-1} =
  R^{\lambda-1}
  \sum_{\satop{\bsu\le\bsm}{|\bsu|\ge\lambda-1}} \bsbeta^\bsu\,\bbB_{|\bsu|,\lambda-1}\,
  \calS(\bsm,\bsu)
\]
for all $\bsm\ne\bszero$ with $|\bsm|\ge \lambda-1$. Recall that
$\bbA_{\bsm,\lambda-1}=0$ for $|\bsm|< \lambda-1$. Then for all $\bsnu\in\indx$ and all $j\ge 1$
we have from \eqref{eq:A1-per} that
\begin{align*}
  \bbA_{\bsnu+\bse_j,\lambda} 
  &= R
  \sum_{\substack{\bsm\le\bsnu \\ |\bsm|\ge\lambda-1}} \binom{\bsnu}{\bsm}
  \bigg(\sum_{\bsw\le\bsnu-\bsm+\bse_j} \bsbeta^{\bsw}\,\Gamma_{|\bsw|}\,
  \calS(\bsnu-\bsm+\bse_j,\bsw)\bigg) \\
  &\qquad\qquad\qquad\times
  R^{\lambda-1}
  \sum_{\substack{\bsu\le\bsm \\ |\bsu|\ge\lambda-1}} \bsbeta^\bsu\,\bbB_{|\bsu|,\lambda-1}\,
  \calS(\bsm,\bsu).
\end{align*}
Next we swap the order of summation twice
\begin{align} \label{eq:sum1}
 \bbA_{\bsnu+\bse_j,\lambda}
 &= R^\lambda
 \sum_{\bsw\le\bsnu+\bse_j} \bsbeta^\bsw\,\Gamma_{|\bsw|}
 \sum_{\substack{\bsm\le\bsnu+\bse_j-\bsw \\ |\bsm|\ge \lambda-1}}
 \binom{\bsnu}{\bsm}\, \calS(\bsnu-\bsm+\bse_j,\bsw) \nonumber\\
 &\qquad\qquad\qquad\times
 \sum_{\substack{\bsu\le\bsm \\ |\bsu|\ge\lambda-1}} \bsbeta^\bsu\,
 \bbB_{|\bsu|,\lambda-1}\, \calS(\bsm,\bsu) \nonumber\\
 &= R^\lambda
 \sum_{\bsw\le\bsnu+\bse_j} \bsbeta^\bsw\,\Gamma_{|\bsw|}
 \sum_{\substack{\bsu\le\bsnu+\bse_j-\bsw \\ |\bsu|\ge\lambda-1}}
 \bsbeta^\bsu\,\bbB_{|\bsu|,\lambda-1} \nonumber\\
 &\qquad\qquad\qquad\times
 \sum_{\substack{\bsu\le\bsm\le\bsnu+\bse_j-\bsw \\ |\bsm|\ge \lambda-1}}
 \binom{\bsnu}{\bsm}\, \calS(\bsnu-\bsm+\bse_j,\bsw)\, \calS(\bsm,\bsu),
\end{align}
where we may drop the condition $|\bsm|\ge \lambda-1$ because it is
already taken care of by the earlier condition $|\bsu|\ge\lambda-1$.

The sum over $\bsm$ in \eqref{eq:sum1} is
\begin{align} \label{eq:sum2}
 \sum_{\bsu\le\bsm\le\bsnu+\bse_j-\bsw} \!\!
 \binom{\bsnu}{\bsm}\, \calS(\bsnu-\bsm+\bse_j,\bsw)\, \calS(\bsm,\bsu)
 = \binom{\bsw+\bsu-\bse_j}{\bsw-\bse_j}\, \calS(\bsnu+\bse_j,\bsw+\bsu),
\end{align}
since (see e.g., \cite[formula~26.8.23]{NIST})
\begin{align*}
 \sum_{m=u}^{\nu-w} \binom{\nu}{m}\, \calS(\nu-m,w)\,\calS(m, u)
 = \binom{w+u}{w}\, \calS(\nu,w+u),
\end{align*}
while
\begin{align*}
 &\sum_{m=u}^{\nu+1-w} \binom{\nu}{m}\, \calS(\nu-m+1,w)\,\calS(m, u) \\
 &= \sum_{m=u}^{\nu+1-w}
 \binom{\nu}{m} \Big(w\, \calS(\nu-m,w) + \calS(\nu-m,w-1)\Big) \,\calS(m, u) \\
 &= w \sum_{m=u}^{\nu-w} \binom{\nu}{m}\, \calS(\nu-m,w) \,\calS(m, u)
 + \sum_{m=u}^{\nu-(w-1)} \binom{\nu}{m}\, \calS(\nu-m,w-1) \,\calS(m, u) \\
 &= w \binom{w+u}{w}\, \calS(\nu,w+u) + \binom{(w-1)+u}{w-1}\, \calS(\nu, w-1+u) \\
 &= \binom{w+u-1}{w-1} \bigg( (w+u)\, \calS(\nu,w+u) + \calS(\nu, w+u-1) \bigg) \\
 &= \binom{w+u-1}{w-1}\, \calS(\nu+1,w+u),
\end{align*}
where we used the recurrence $\calS(n+1,k) = k\,\calS(n,k) + \calS(n,k-1)$ for $0<k<n$ in the first 
and last equalities.
Substituting \eqref{eq:sum2} into \eqref{eq:sum1}, we arrive at
\begin{align*}
 &\bbA_{\bsnu+\bse_j,\lambda}
 = R^\lambda
 \sum_{\bsw\le\bsnu+\bse_j} \bsbeta^\bsw\,\Gamma_{|\bsw|} \!\!
 \sum_{\substack{\bsu\le\bsnu+\bse_j-\bsw \\ |\bsu|\ge\lambda-1}} \!\!
 \bsbeta^\bsu\,\bbB_{|\bsu|,\lambda-1}
 \binom{\bsw+\bsu-\bse_j}{\bsw-\bse_j}\, \calS(\bsnu+\bse_j,\bsw+\bsu) \\
 &= R^\lambda
 \sum_{\bsw\le\bsnu+\bse_j} \bsbeta^\bsw\,\Gamma_{|\bsw|}
 \sum_{\substack{\bsw\le\bsu'\le\bsnu+\bse_j \\ |\bsu'|\ge|\bsw|+\lambda-1}} \!\!
 \beta^{\bsu'-\bsw}\,\bbB_{|\bsu'-\bsw|,\lambda-1}\,
 \binom{\bsu'-\bse_j}{\bsw-\bse_j}\, \calS(\bsnu+\bse_j,\bsu') \\
 &= R^\lambda\,
 \sum_{\satop{\bsu\le\bsnu+\bse_j}{|\bsu|\ge\lambda}} \bsbeta^\bsu
 \bigg(\sum_{\substack{\bsw\le\bsu \\ |\bsw|\le|\bsu|-(\lambda-1)}}
 \binom{\bsu-\bse_j}{\bsw-\bse_j}\, \Gamma_{|\bsw|}\,
 \bbB_{|\bsu-\bsw|,\lambda-1}\bigg)\,\calS(\bsnu+\bse_j,\bsu),
\end{align*}
where we substituted $\bsu' = \bsw+\bsu$ and then swapped the order of
summation (and relabeled $\bsu'$ back to $\bsu$ at the same time). We see
that implicitly $\bsw\ge\bse_j$ (otherwise the binomial coefficient will
be $0$), which leads to the condition $|\bsu|\ge\lambda$. We can rewrite
the inner sum over $\bsw$ as
\begin{align*}
 \sum_{\substack{\bsw\le\bsu \\ |\bsw|\le|\bsu|-(\lambda-1)}} \!\!\!\!\!
 \binom{\bsu-\bse_j}{\bsw-\bse_j}\, \Gamma_{|\bsw|}\,
 \bbB_{|\bsu-\bsw|,\lambda-1}
 &=
 \sum_{i=1}^{|\bsu|-(\lambda-1)}
 \sum_{\substack{\bse_j\le\bsw\le\bsu \\ |\bsw|=i}}
 \binom{\bsu-\bse_j}{\bsw-\bse_j}\, \Gamma_i\,\bbB_{|\bsu|-i,\lambda-1}\\
 &=
 \sum_{i=1}^{|\bsu|-(\lambda-1)}
 \binom{|\bsu|-1}{i-1}\, \Gamma_{i}\,\bbB_{|\bsu|-i,\lambda-1}\\
 &=
 \sum_{i'=\lambda-1}^{|\bsu| - 1} \binom{|\bsu|-1}{i'}\, \Gamma_{|\bsu|-i'}\,
 \bbB_{i',\lambda-1}
 = \bbB_{|\bsu|,\lambda},
\end{align*}
where we substituted $i' = |\bsu| - i$ and used the symmetry of the binomial coefficient and the recursive definition of $\bbB_{|\bsu|,\lambda}$ in
\eqref{eq:Bn}. In the second equality above we used \eqref{eq:comb} to obtain
\[
 \sum_{\substack{\bse_j\le\bsw\le\bsu \\ |\bsw|=i}}
 \binom{\bsu-\bse_j}{\bsw-\bse_j}
 = \sum_{\substack{\bsw'\le\bsu-\bse_j \\ |\bsw'|=i-1}}
 \binom{\bsu-\bse_j}{\bsw'}
 = \binom{|\bsu|-1}{i-1}.
\]
Hence we conclude that
\begin{align*}
 \bbA_{\bsnu+\bse_j,\lambda}
 &= R^\lambda
 \sum_{\satop{\bsu\le\bsnu+\bse_j}{|\bsu|\ge\lambda}} \bsbeta^\bsu\,
 \bbB_{|\bsu|,\lambda}\,
 \calS(\bsnu+\bse_j,\bsu),
\end{align*}
which proves the induction step for \eqref{eq:A2-per}. This completes the
induction proof.

If the equality in \eqref{eq:A1-per} is replaced by an inequality $\le$\,, then it is easy
to see from the proof that this will lead to an inequality $\le$ in the final
result \eqref{eq:A2-per}.
\end{proof}

\begin{proofof}{Theorem~\ref{thm:der}(b)}
The proof for the periodic case is more complicated, hence it is helpful as a
preparation step to first consider the case with no hidden layer
\begin{align*}
  G_\theta^{[0]}(\bsy)
  = W_0\, \sin(2\pi\bsy) + \bsv_0
  = \bigg[ \sum_{q=1}^{d_0} W_{0,p,q}\, \sin(2\pi y_q) + v_{0,p}\bigg]_{p=1}^{d_1},
\end{align*}
where $d_0 = s$. For $\bsnu\ne\bszero$ and any index $p$, we may
differentiate with respect to a single variable $y_j$ multiple times but
all cross derivatives vanish:
\begin{align} \label{eq:der0}
  \partial^\bsnu G_\theta^{[0]}(\bsy)_p
  &=
  \begin{cases}
  W_{0,p,j}\, (2\pi)^k \sin(2\pi y_j + \tfrac{k\pi}{2}) & \mbox{if } \bsnu = k \bse_j,\, k\ge 1, \\
  0 & \mbox{otherwise},
  \end{cases}
\end{align}
from which we conclude that
\begin{align} \label{eq:bnd0}
  |\partial^{\bsnu} G_\theta^{[0]}(\bsy)_{p}|
  &\le
  \begin{cases}
  \beta_j\, (2\pi)^{k} & \mbox{if } \bsnu = k \bse_j,\, k\ge 0, \\
  0 & \mbox{otherwise},
  \end{cases}
\end{align}
where we used the assumption \eqref{eq:beta}.

We are ready to consider a shallow network with one hidden layer
\[
  G_\theta^{[1]}(\bsy) = W_1\,\sigma (G_\theta^{[0]}(\bsy)) + \bsv_1
  = \bigg[
    \sum_{q=1}^{d_1} W_{1,p,q}\, \sigma (G_\theta^{[0]}(\bsy)_q ) + v_{1,p}
    \bigg]_{p=1}^{d_2}.
\]
For any $\bsnu\ne\bszero$ and any index $p$, we have
\begin{align} \label{eq:der7a}
  \partial^\bsnu G_\theta^{[1]}(\bsy)_p
  = \sum_{q=1}^{d_1} W_{1,p,q}\,
  \partial^\bsnu \big( \sigma (G_\theta^{[0]}(\bsy)_q)\big).
\end{align}
For each index $q$, the recursive form of the multivariate Fa\`a di Bruno
formula gives (see \cite[formula (3.1) and (3.5)]{Savits06} with the
dimensionality $m$ there set to $1$)
\begin{align} \label{eq:der7b}
  \partial^\bsnu \big( \sigma (G_\theta^{[0]}(\bsy)_q)\big)
  &=
  \sum_{\lambda=1}^{|\bsnu|} \sigma^{(\lambda)} (G_\theta^{[0]}(\bsy)_q)\,
  \alpha_{\bsnu,\lambda}(\bsy),
\end{align}
where the auxiliary functions $\alpha_{\bsnu,\lambda}(\bsy)$ are
defined recursively by $\alpha_{\bsnu,0}(\bsy) :=
\delta_{\bsnu,\bszero}$, $\alpha_{\bsnu,\lambda}(\bsy) := 0$ for
$\lambda> |\bsnu|$, and otherwise
\[
  \alpha_{\bsnu+\bse_j,\lambda}(\bsy) :=
  \sum_{\bsm\le\bsnu} \binom{\bsnu}{\bsm}
  \partial^{\bsm+\bse_j} G_\theta^{[0]}(\bsy)_q \,
  \alpha_{\bsnu-\bsm,\lambda-1}(\bsy).
\]
(Note that we have written the recursion differently from the non-periodic
case, with $\bsm$ and $\bsnu-\bsm$ swapped over, for convenience here.)
Using \eqref{eq:bnd0}, we obtain
\begin{align*}
  |\alpha_{\bsnu+\bse_j,\lambda}(\bsy)|
  &\le
  \sum_{\bsm\le\bsnu} \binom{\bsnu}{\bsm}
  |\partial^{\bsm+\bse_j} G_\theta^{[0]}(\bsy)_{q}| \, |\alpha_{\bsnu-\bsm,\lambda-1}(\bsy)| \\
  &\le \beta_j \sum_{k=0}^{\nu_j} (2\pi)^{k+1}\, \binom{\nu_j}{k}\,
 |\alpha_{\bsnu-k\bse_j,\lambda-1}(\bsy)|.
\end{align*}
To solve this recurrence, we apply
\cite[Lemma~A.4]{HHKKS24} (taking there $c = 2\pi$ and $d=1$
which is the dimensionality of $\bslam$ there) to obtain the bound
\begin{align} \label{eq:alpha-per0}
  |\alpha_{\bsnu,\lambda}(\bsy)|
  \le (2\pi)^{|\bsnu|}
  \sum_{\satop{\bsm\le\bsnu}{|\bsm|=\lambda}} \bsbeta^\bsm\,\calS(\bsnu,\bsm).
\end{align}
(This bound is tight in the sense that it would be an equality if the recurrence were an equality.)
Applying \eqref{eq:sigma} and \eqref{eq:alpha-per0} in \eqref{eq:der7b}, we obtain for all $\bsnu\ne\bszero$,
\begin{align} \label{eq:der7c}
  &\big|\partial^\bsnu \big( \sigma (G_\theta^{[0]}(\bsy)_q)\big)\big|
  \le \sum_{\lambda=1}^{|\bsnu|} A_\lambda\, (2\pi)^{|\bsnu|}
  \sum_{\satop{\bsm\le\bsnu}{|\bsm|=\lambda}} \bsbeta^\bsm\, \calS(\bsnu,\bsm) \nonumber\\
  &= (2\pi)^{|\bsnu|}
  \sum_{\bszero\ne\bsm\le\bsnu} \bsbeta^\bsm\, A_{|\bsm|}\, \calS(\bsnu,\bsm)
  = (2\pi)^{|\bsnu|}
  \sum_{\bsm\le\bsnu} \bsbeta^\bsm\, \Gamma_{|\bsm|}^{[1]}\, \calS(\bsnu,\bsm).
\end{align}
The $\bsm=\bszero$ term is $0$ and so it was added back into the sum. (This is
because $\bsnu\ne\bszero$ means there exists a $\nu_j >0$ and $\calS(\nu_j,0) =
0$, and hence $\calS(\bsnu,\bszero)=0$.) We also used $A_{|\bsm|} = \Gamma_{|\bsm|}^{[1]}$ from
\eqref{eq:G-def}. 
Applying \eqref{eq:R} and \eqref{eq:der7c} in \eqref{eq:der7a}, we obtain for all $\bsnu\ne\bszero$,
\begin{align*}
  |\partial^\bsnu G_\theta^{[1]}(\bsy)_p|
  \le (2\pi)^{|\bsnu|}\,R_1\,
  \sum_{\bsm\le\bsnu} \bsbeta^\bsm\, \Gamma_{|\bsm|}^{[1]}\, \calS(\bsnu,\bsm).
\end{align*}
Hence \eqref{eq:der-per} holds for $\ell=1$.

Now consider $\ell \ge 2$ and suppose the result \eqref{eq:der-per} holds
for depth $\ell-1$. We can write
\begin{align*}
  G_\theta^{[\ell]}(\bsy) = W_\ell\, \sigma (G_\theta^{[\ell-1]}(\bsy)) + \bsv_\ell
  = \bigg[ \sum_{q=1}^{d_\ell} W_{\ell,p,q}\,
  \sigma (G_\theta^{[\ell-1]}(\bsy)_q) + v_{\ell,p} \bigg]_{p=1}^{d_{\ell+1}}.
\end{align*}
For any $\bsnu\ne\bszero$ and any index $p$, we have
\begin{align} \label{eq:der8a}
  \partial^\bsnu G_\theta^{[\ell]}(\bsy)_p
  &= \sum_{q=1}^{d_\ell} W_{\ell,p,q}\,
 \partial^\bsnu \big( \sigma (G_\theta^{[\ell-1]}(\bsy)_q)\big).
\end{align}
For each index $q$, the recursive form of the multivariate Fa\`a di Bruno
formula gives (see \cite[formula (3.1) and (3.5)]{Savits06} with the
dimensionality $m$ there set to $1$)
\begin{align} \label{eq:der8b}
  \partial^\bsnu \big( \sigma (G_\theta^{[\ell-1]}(\bsy)_q)\big)
  &=
  \sum_{\lambda=1}^{|\bsnu|} \sigma^{(\lambda)} (G_\theta^{[\ell-1]}(\bsy)_q)\,
  \alpha_{\bsnu,\lambda}(\bsy),
\end{align}
where now the auxiliary functions $\alpha_{\bsnu,\lambda}(\bsy)$
are defined recursively by $\alpha_{\bsnu,0}(\bsy) :=
\delta_{\bsnu,\bszero}$, $\alpha_{\bsnu,\lambda}(\bsy) := 0$
for $\lambda
> |\bsnu|$, and otherwise
\[
  \alpha_{\bsnu+\bse_j,\lambda}(\bsy) :=
  \sum_{\bsm\le\bsnu} \binom{\bsnu}{\bsm}
  \partial^{\bsnu-\bsm+\bse_j} G_\theta^{[\ell-1]}(\bsy)_q \,
  \alpha_{\bsm,\lambda-1}(\bsy).
\]
Using the induction hypothesis \eqref{eq:der-per} for depth $\ell-1$, we
obtain
\begin{align} \label{eq:tri-per}
  &|\alpha_{\bsnu+\bse_j,\lambda}(\bsy)|
  \le \sum_{\bsm\le\bsnu} \binom{\bsnu}{\bsm}
  |\partial^{\bsnu-\bsm+\bse_j} G_\theta^{[\ell-1]}(\bsy)_q| \,
  |\alpha_{\bsm,\lambda-1}(\bsy)| \\
  &\le \sum_{\bsm\le\bsnu } \binom{\bsnu}{\bsm}
  \bigg(
  (2\pi)^{|\bsnu-\bsm|+1}\,R_{\ell-1} \!\!\!\!\!
  \sum_{\bsw\le\bsnu-\bsm+\bse_j} \!\!\!\!
  \bsbeta^{\bsw}\,\Gamma_{|\bsw|}^{[\ell-1]}\,
  \calS(\bsnu-\bsm+\bse_j,\bsw)\bigg)
  |\alpha_{\bsm,\lambda-1}(\bsy)|. \nonumber
\end{align}
Taking $R = R_{\ell-1}$ and $\Gamma_n = \Gamma_n^{[\ell-1]}$ in Lemma~\ref{lem:recur-per} (with inequalities $\le$) then yields for all
$\bsnu\ne\bszero$,
\begin{align} \label{eq:alpha-per}
  |\alpha_{\bsnu,\lambda}(\bsy)|
  \le R_{\ell-1}^\lambda\,(2\pi)^{|\bsnu|}
  \sum_{\satop{\bsu\le\bsnu}{|\bsu|\ge\lambda}}
  \bsbeta^{\bsu}\,\bbB_{|\bsu|,\lambda}^{[\ell-1]}\, \calS(\bsnu,\bsu),
\end{align}
with the sequence $\bbB_{n,\lambda}^{[\ell-1]}$ as defined in \eqref{eq:B-def}. 
Applying \eqref{eq:sigma} and \eqref{eq:alpha-per} in \eqref{eq:der8b}, we obtain
\begin{align} \label{eq:der8c}
  \big|\partial^\bsnu \big( \sigma (G_\theta^{[\ell-1]}(\bsy)_q)\big)\big|
  &\le \sum_{\lambda=1}^{|\bsnu|} A_\lambda\,
  R_{\ell-1}^\lambda\,(2\pi)^{|\bsnu|}\,
  \sum_{\satop{\bsu\le\bsnu}{|\bsu|\ge\lambda}}
  \bsbeta^{\bsu}\,\bbB_{|\bsu|,\lambda}^{[\ell-1]} \, \calS(\bsnu,\bsu) \nonumber\\
  &= (2\pi)^{|\bsnu|}\, \sum_{\bsu\le\bsnu} \bsbeta^{\bsu}
  \bigg(
  \underbrace{
  \sum_{\lambda=1}^{|\bsu|} A_\lambda\, R_{\ell-1}^\lambda\,\bbB_{|\bsu|,\lambda}^{[\ell-1]}
  }_{=\,\Gamma_{|\bsu|}^{[\ell]}}
  \bigg)\, \calS(\bsnu,\bsu),
\end{align}
with the sequence $\Gamma_n^{[\ell]}$ as defined in \eqref{eq:B-def}. Applying \eqref{eq:R} and \eqref{eq:der8c} in \eqref{eq:der8a}, we arrive at
\begin{align*}
  |\partial^\bsnu G_\theta^{[\ell]}(\bsy)_p|
  \le (2\pi)^{|\bsnu|}\, R_\ell\,
  \sum_{\bsu\le\bsnu} \bsbeta^{\bsu}\,\Gamma_{|\bsu|}^{[\ell]}\, \calS(\bsnu,\bsu).
\end{align*}
This completes the proof.
\end{proofof}

\subsection{Proofs for popular activation functions}
\label{sec:act}

The \textbf{logistic sigmoid} function
\begin{align} \label{eq:sigmoid}
  \sigma(x) = \frac{1}{1+e^{-x}}
\end{align}
is monotone on $\bbR$ and has range $(0,1)$. From \cite[formula
(15)]{MinWil93} we have for $n\ge 1$
\begin{align} \label{eq:sigmoid-n}
  \sigma^{(n)}(x) = \sum_{k=1}^n (-1)^{k-1} E(n,k-1)\,[\sigma(x)]^k\,[1-\sigma(x)]^{n+1-k},
\end{align}
where $E(n,k) := \sum_{i=0}^k (-1)^i \binom{n+1}{i} (k+1-i)^n$ is the
Eulerian number, i.e., the number of permutations of the integers $1$ to
$n$ in which exactly $k$ elements are greater than the previous element.
Thus
\[
  |\sigma^{(n)}(x)| \le \sum_{k=1}^n E(n,k-1) = \sum_{k=0}^{n-1} E(n,k) = n!\,,
\]
where conveniently the sum of Eulerian numbers is equal to $n!$\,.

We can actually obtain a tighter upper bound. Indeed, the function $f_k(y)
:= y^k(1-y)^{n+1-k}$ for $1\le k\le n$ vanishes when $y=0$ and $y=1$; it
is positive on the interval $(0,1)$ with a unique maximum at $y =
\frac{k}{n+1}$. Furthermore, $f_k(\frac{k}{n+1})$ is largest when $k=1$ or
$n$, and this value is $f_1(\frac{1}{n+1}) = \frac{n^n}{(n+1)^{n+1}}$. Thus
for all $n\ge 1$ we have
\begin{align} \label{eq:sn-bound}
  |\sigma^{(n)}(x)| \le f_1\Big(\frac{1}{n+1}\Big) \sum_{k=1}^n E(n,k-1)
  = \frac{n^n}{(n+1)^{n+1}}\sum_{k=0}^{n-1} E(n,k)
  &= \frac{n^n}{(n+1)^{n+1}}\,n! \,.
\end{align}

The \textbf{hyperbolic tangent} activation function 
is related to the logistic sigmoid \eqref{eq:sigmoid} by
\[
  \tanh(x) = \frac{e^x-e^{-x}}{e^x+e^{-x}} = 2\,\sigma(2x)-1,
\]
and therefore for $n\ge1$ we have
\[
 \frac{\rd^n}{\rd x^n}\tanh(x) = 2^{n+1}\,\sigma^{(n)}(2x),
\]
which combined with \eqref{eq:sn-bound} gives
\[
  \Big|\frac{\rd^n}{\rd x^n}\tanh(x)\Big| \le 2^{n+1}\,
  \frac{1}{n} \Big(\frac{n}{n+1}\Big)^{n+1} n!
  < 2^n\,n!\qquad\mbox{for all $n\ge 1$.}
\]

The \textbf{swish} activation function is related to the logistic sigmoid function \eqref{eq:sigmoid}~by
\[
  {\rm swish}(x) = \frac{x}{1+e^{-x}} = x\,\sigma(x),
\]
and therefore
\[
 \frac{\rd^n}{\rd x^n}{\rm swish}(x)
 = x\,\sigma^{(n)}(x) + n\,\sigma^{(n-1)}(x).
\]
Using \eqref{eq:sigmoid-n} we obtain
\begin{align*}
 x\,\sigma^{(n)}(x)
 &= \sum_{k=1}^n (-1)^{k-1} E(n,k-1)\,x\, [\sigma(x)]^k\,[1-\sigma(x)]^{n+1-k} \\
 &= \sum_{k=1}^n (-1)^{k-1} E(n,k-1)\,g(x)\, [\sigma(x)]^{k-1}\,[1-\sigma(x)]^{n-k},
\end{align*}
where (noting that $k-1\ge 0$ and $n-k\ge 0$) we defined
\[
 g(x) := x\,\sigma(x)\,(1-\sigma(x))
 = \frac{x\, e^{-x}}{(1+e^{-x})^2}
 = \frac{x}{(e^{x/2}+e^{-x/2})^2}.
\]
Clearly $g$ is an odd function and $g'(x) =
e^{-x}(1-x+e^{-x}+x\,e^{-x})/(1+e^{-x})^3$ vanishes exactly when $x\approx
\pm 1.543$, leading to two global extrema, and therefore $|g(x)|\le
g_{\max} \approx 0.224$. Thus
\begin{align} \label{eq:xsn}
 |x\,\sigma^{(n)}(x)|
 \le g_{\max}\sum_{k=1}^n E(n,k-1)
 = g_{\max}\sum_{k=0}^{n-1} E(n,k) = g_{\max}\,n!\,.
\end{align}
Hence we can use \eqref{eq:xsn} and \eqref{eq:sn-bound} to obtain for
$n\ge 2$,
\begin{align*}
 \Big|\frac{\rd^n}{\rd x^n}{\rm swish}(x)\Big|
 &\le |x\,\sigma^{(n)}(x)| + n\, |\sigma^{(n-1)}(x)| \\
 &\le g_{\max}\,n! + n\,\frac{(n-1)^{n-1}}{n^n} (n-1)! \\
 &= g_{\max}\,n! + \frac{1}{n}\Big(\frac{n-1}{n}\Big)^{n-1}\, n!
 < \Big(\frac{1}{4} + \frac{1}{n} \Big) n!
 < n!\,.
\end{align*}

For $n=1$ we have
\[
 \frac{\rd}{\rd x}{\rm swish}(x)
 = x\,\sigma'(x) + \sigma(x)
 = x\,\frac{e^{-x}}{(1+e^{-x})^2} + \frac{1}{1+e^{-x}}
 = \frac{1+e^{-x}+x\,e^{-x}}{(1+e^{-x})^2},
\]
and
\begin{align*}
 \frac{\rd^2}{\rd x^2}{\rm swish}(x)
 &= \frac{-x\,e^{-x}}{(1+e^{-x})^2} +
 \frac{2\,e^{-x} (1+e^{-x}+x\,e^{-x})}{(1+e^{-x})^3} \\
 &= \frac{e^{-x}[-x(1+e^{-x}) + 2(1+e^{-x}+x\,e^{-x})]}{(1+e^{-x})^3} \\
 &= \frac{e^{-x}[-x + 2 + 2e^{-x}+ x\,e^{-x}]}{(1+e^{-x})^3};
\end{align*}
the latter vanishes at $x\approx \pm 2.3994$, giving local maximum/mininum
values $1.0998$ and $-0.0998$ for the first derivative. Thus
$|\frac{\rd}{\rd x}{\rm swish}(x)| \le 1.1$. Combining this with the upper
bound $n!$ for $n\ge 2$, we conclude that
\begin{align*}
 \Big|\frac{\rd^n}{\rd x^n}{\rm swish}(x)\Big|
 \le 1.1\,n! \qquad\mbox{for all $n\ge 1$.}
\end{align*}

In summary, we have shown that the bound $A_n$ in \eqref{eq:sigma} can be taken as
\begin{align} \label{eq:An-cases}
  A_n =
  \begin{cases}
  n! & \text{for sigmoid}, \\
  2^n\,n! & \text{for tanh}, \\
  1.1\,n! & \text{for swish},
  \end{cases}
\end{align}
leading us to consider the general form $A_n = \xi\,\tau^n\,n!$\,.

\begin{lemma} \label{lem:complex}
Let $\xi>0$ and $\tau>0$. If $A_n = \xi\,\tau^n\,n!$ for
$n\ge 1$ and $R_\ell > 0$ for $\ell\ge 1$, then the sequences defined
recursively in \eqref{eq:G-def}--\eqref{eq:B-def} are given explicitly by
\begin{align}
  \Gamma_n^{[\ell]}
  = \bbB_{n,1}^{[\ell]}
  &= P_{\ell-1}\,\bigg(\sum_{k=0}^{\ell-1} P_k\bigg)^{n-1}\,\xi\,\tau^n\,n!
  \quad\mbox{for $\ell\ge 1$ and $n\ge 1$}, \label{eq:G-exp} \\
  \bbB_{n,\lambda}^{[\ell]}
  &= P_{\ell-1}^\lambda\,\bigg(\sum_{k=0}^{\ell-1} P_k\bigg)^{n-\lambda}\,
  \xi^\lambda\,\tau^n\,\frac{n!}{\lambda!} \binom{n-1}{\lambda-1}
  \quad\mbox{for $\ell\ge 1$ and $n\ge\lambda\ge 1$}, \label{eq:B-exp}
\end{align}
where
\begin{align} \label{eq:P-def}
  P_0 := 1, \qquad
  P_k := \prod_{t=1}^k (\xi\,\tau R_t) \qquad\mbox{for $k\ge 1$}.
\end{align}
\end{lemma}

\begin{proof}
For $\ell=1$, we have from \eqref{eq:G-def}--\eqref{eq:B-def} that
$\Gamma_n^{[1]} = \bbB_{n,1}^{[1]} = A_n = \xi\,\tau^n\,n!$ for all $n\ge
1$, which agrees with \eqref{eq:G-exp} since $P_0 = 1$ from
\eqref{eq:P-def}. Thus the formula \eqref{eq:G-exp} holds for $\ell=1$ and
all $n\ge 1$. And trivially, the formula \eqref{eq:B-exp} holds for
$\ell=1$ and $\lambda=1$ and all $n\ge\lambda = 1$.

For $\ell=1$ and some $\lambda\ge 2$, assume that the formula
\eqref{eq:B-exp} for $\bbB_{n,\lambda'}^{[1]}$ holds for all $1\le
\lambda' < \lambda$ and all $n\ge \lambda'$. Then from the recursion
\eqref{eq:B-def} we have
\begin{align*}
 \bbB_{n,\lambda}^{[1]}
 &= \sum_{i=\lambda-1}^{n-1} \binom{n-1}{i}\,
 \Gamma_{n-i}^{[1]}\,\bbB_{i,\lambda-1}^{[1]} \\
 &= \sum_{i=\lambda-1}^{n-1} \binom{n-1}{i}\,
 \xi\,\tau^{n-i}\,(n-i)!\,\xi^{\lambda-1}\,\tau^i\,
 \frac{i!}{(\lambda-1)!} \binom{i-1}{\lambda-2} \\
 &= \xi^\lambda\,\tau^n\,\frac{(n-1)!}{(\lambda-1)!}\,
 \sum_{i=\lambda-1}^{n-1} (n-i) \binom{i-1}{\lambda-2},
\end{align*}
with
\begin{align} \label{eq:hockey}
 \sum_{i=\lambda-1}^{n-1} (n-i)\,\binom{i-1}{\lambda-2}
 &= n \sum_{i=\lambda-1}^{n-1} \binom{i-1}{\lambda-2}
  - (\lambda-1) \sum_{i=\lambda-1}^{n-1} \binom{i}{\lambda-1} \nonumber\\
 &= n \sum_{i=\lambda-2}^{n-2} \binom{i}{\lambda-2}
  - (\lambda-1) \sum_{i=\lambda-1}^{n-1} \binom{i}{\lambda-1} \nonumber\\
 &= n \binom{n-1}{\lambda-1} - (\lambda-1) \binom{n}{\lambda}
 = \lambda \binom{n}{\lambda} - (\lambda-1) \binom{n}{\lambda}
 = \binom{n}{\lambda},
\end{align}
where we used the ``hockey-stick identity'' $\sum_{m=k}^n \binom{m}{k} =
\binom{n+1}{k+1}$. Thus we have
\begin{align} \label{eq:same}
 \bbB_{n,\lambda}^{[1]}
 = \xi^\lambda\,\tau^n\,\frac{(n-1)!}{(\lambda-1)!}\, \binom{n}{\lambda}
 = \xi^\lambda\,\tau^n\,\frac{n!}{\lambda!}\, \binom{n-1}{\lambda-1}.
\end{align}
Hence by induction on $\lambda$, the formula \eqref{eq:B-exp} holds for
$\ell=1$ and all $n\ge\lambda\ge 1$.

Consider now a fixed level $\ell\ge 2$ and assume that the formulas
\eqref{eq:G-exp}--\eqref{eq:B-exp} hold in level $\ell-1$. Then for all
$n\ge 1$ we have from \eqref{eq:G-def} that
\begin{align*}
 \Gamma_n^{[\ell]}
 = \bbB_{n,1}^{[\ell]}
 &= \sum_{\lambda=1}^n A_\lambda\, R_{\ell-1}^\lambda\, \bbB_{n,\lambda}^{[\ell-1]} \\
 &= \sum_{\lambda=1}^n \xi\,\tau^\lambda\,\lambda!\, R_{\ell-1}^\lambda\,
 P_{\ell-2}^\lambda\,\bigg(\sum_{k=0}^{\ell-2} P_k\bigg)^{n-\lambda}\,
 \xi^\lambda\,\tau^n\,
 \frac{n!}{\lambda!}\,\binom{n-1}{\lambda-1} \\
 &= \xi\,\tau^n\,n!\, \sum_{\lambda=1}^n P_{\ell-1}^\lambda\,
 \bigg(\sum_{k=0}^{\ell-2} P_k\bigg)^{n-\lambda}\, \binom{n-1}{\lambda-1} \\
 &= \xi\,\tau^n\,n!\, P_{\ell-1} \sum_{\lambda'=0}^{n-1} P_{\ell-1}^{\lambda'}\,
 \bigg(\sum_{k=0}^{\ell-2} P_k\bigg)^{n-1-\lambda'}\, \binom{n-1}{\lambda'} \\
 &= \xi\,\tau^n\,n!\, P_{\ell-1}\,
 \bigg(P_{\ell-1} + \sum_{k=0}^{\ell-2} P_k\bigg)^{n-1}
 = P_{\ell-1} \bigg(\sum_{k=0}^{\ell-1} P_k\bigg)^{n-1}\,\xi\,\tau^n\,n!\,,
\end{align*}
where in the fourth equality we used $(\xi\,\tau\,R_{\ell-1})\,P_{\ell-2}
= P_{\ell-1}$ from \eqref{eq:P-def}, and in the sixth equality we used the
binomial theorem. Thus the formula \eqref{eq:B-exp} holds for this level
$\ell\ge 2$ and all $n\ge 1$. Also trivially, the formula \eqref{eq:B-exp}
holds for this level $\ell\ge 2$ and $\lambda=1$, and all $n\ge\lambda$, because \eqref{eq:B-exp} coincides with \eqref{eq:G-exp} when $\lambda = 1$.

For the same level $\ell\ge 2$ and some $\lambda\ge 2$, assume that the
formula \eqref{eq:B-exp} for $\bbB_{n,\lambda'}^{[\ell]}$ holds for all
$1\le\lambda'< \lambda$ and all $n\ge\lambda'$. Then from the recursion
\eqref{eq:B-def} we have
\begin{align*}
 \bbB_{n,\lambda}^{[\ell]}
 &= \sum_{i=\lambda-1}^{n-1} \binom{n-1}{i}\,
 \Gamma_{n-i}^{[\ell]}\,\bbB_{i,\lambda-1}^{[\ell]} \\
 &= \sum_{i=\lambda-1}^{n-1} \binom{n-1}{i}\,
 P_{\ell-1}\,\bigg(\sum_{k=0}^{\ell-1} P_k\bigg)^{n-i-1}\,\xi\,\tau^{n-i}\,(n-i)!\\
 &\qquad\qquad\qquad\times
 P_{\ell-1}^{\lambda-1}\,\bigg(\sum_{k=0}^{\ell-1} P_k\bigg)^{i-(\lambda-1)}\,
 \xi^{\lambda-1}\,\tau^i\,\frac{i!}{(\lambda-1)!} \binom{i-1}{\lambda-2} \\
 &= P_{\ell-1}^{\lambda}\,\bigg(\sum_{k=0}^{\ell-1} P_k\bigg)^{n-\lambda}\,
 \xi^\lambda\,\tau^n\,\frac{(n-1)!}{(\lambda-1)!}
 \sum_{i=\lambda-1}^{n-1} (n-i)\,\binom{i-1}{\lambda-2} \\
 &= P_{\ell-1}^\lambda\,\bigg(\sum_{k=0}^{\ell-1} P_k\bigg)^{n-\lambda}\,
 \xi^\lambda\,\tau^n\,
 \frac{n!}{\lambda!}\,\binom{n-1}{\lambda-1},
\end{align*}
where the last step follows from \eqref{eq:hockey} and the second equality in 
\eqref{eq:same}. This is exactly \eqref{eq:B-exp}, thus by induction on $\lambda$, the formula \eqref{eq:B-exp} holds
for this level $\ell\ge 2$ and all $n\ge\lambda\ge 1$.

Hence by induction on the level $\ell$, we conclude that the formulas
\eqref{eq:G-exp}--\eqref{eq:B-exp} hold for all levels $\ell\ge 1$ and all
$n\ge\lambda\ge 1$. This completes the proof.
\end{proof}

\subsection{Proofs for bounds on norms} \label{sec:norm}

\begin{proofof}{Theorem~\ref{thm:diff}} Consider first the
non-periodic case. With the restrictions
\eqref{eq:demand1}--\eqref{eq:demand3}, both the target function $G(\bsy)$
and the DNN $G^{[L]}_\theta(\bsy)$ satisfy the regularity bound
\eqref{eq:tar-np}, and thus their difference satisfies
\[
  \big|\partial^{\bsnu} \big(G(\bsy)_p - G^{[L]}_\theta(\bsy)_p\big)\big|
  \le \big|\partial^{\bsnu} G(\bsy)_p \big|
  + \big|\partial^{\bsnu} G^{[L]}_\theta(\bsy)_p \big|
  \le 2\,C\,|\bsnu|!\,\bsb^\bsnu.
\]
For the square of the difference, we use the Leibniz product rule to
obtain
\begin{align*}
 &\big|\partial^{\bsnu} \big(G(\bsy)_p - G^{[L]}_\theta(\bsy)_p\big)^2\big| \\
 &\le \sum_{\bsm\le\bsnu} \binom{\bsnu}{\bsm}
 \big|\partial^{\bsm} \big(G(\bsy)_p - G^{[L]}_\theta(\bsy)_p\big)\big| \,
 \big|\partial^{\bsnu - \bsm} \big(G(\bsy)_p - G^{[L]}_\theta(\bsy)_p\big)\big| \\
 &\le \sum_{\bsm\le\bsnu} \binom{\bsnu}{\bsm}
 \big(2\,C\, |\bsm|!\, \bsb^\bsm\big)
 \big(2\,C\, |\bsnu-\bsm|!\, \bsb^{\bsnu-\bsm} \big)\\
 &= 4\,C^2\,(|\bsnu|+1)!\, \bsb^\bsnu\,,
\end{align*}
where we used the identity (see e.g., \cite[formula~(9.4)]{KN16})
\begin{align} \label{eq:1sum}
 \sum_{\bsm\le\bsnu} \binom{\bsnu}{\bsm}\,|\bsm|!\,|\bsnu-\bsm|!
 = (|\bsnu|+1)!\,.
\end{align}
This yields the norm bound \eqref{eq:final-np} as required.

Analogously, for the periodic case both the target function $G(\bsy)$ and
the DNN $G^{[L]}_\theta(\bsy)$ satisfy the regularity bound
\eqref{eq:tar-per}. Thus we have
\begin{align*}
  &\big|\partial^{\bsnu} \big(G(\bsy)_p - G^{[L]}_\theta(\bsy)_p\big)\big| \\
  &\le \big|\partial^{\bsnu} G(\bsy)_p \big|
  + \big|\partial^{\bsnu} G^{[L]}_\theta(\bsy)_p \big|
  \le 2\,C\,(2\pi)^{|\bsnu|} \sum_{\bsm\le\bsnu} |\bsm|!\,\bsb^\bsm\,\calS(\bsnu,\bsm),
\end{align*}
and
\begin{align*}
 & \big|\partial^{\bsnu} \big(G(\bsy)_p - G^{[L]}_\theta(\bsy)_p\big)^2\big| \\
 &\le \sum_{\bsm\le\bsnu} \binom{\bsnu}{\bsm}
 \big|\partial^{\bsm} \big(G(\bsy)_p - G^{[L]}_\theta(\bsy)_p\big)\big| \,
 \big|\partial^{\bsnu - \bsm} \big(G(\bsy)_p - G^{[L]}_\theta(\bsy)_p\big)\big|
 \\
 &\le \sum_{\bsm\le\bsnu} \binom{\bsnu}{\bsm}
 \bigg(2\,C\,(2\pi)^{|\bsm|} \sum_{\bsw\le\bsm} |\bsw|!\,\bsb^\bsw\,\calS(\bsm,\bsw)\!\bigg) \\
 &\qquad\qquad\times
 \bigg(2\,C\,(2\pi)^{|\bsnu-\bsm|} \sum_{\bsu\le\bsnu-\bsm}
 |\bsu|!\,\bsb^\bsu\, \calS(\bsnu-\bsm,\bsu)\!\bigg) \\
 &= 4\,C^2\,(2\pi)^{|\bsnu|} \sum_{\bsm\le\bsnu}
 \bigg(\sum_{\bsw\le\bsm} \binom{\bsm}{\bsw}\,|\bsw|!\,|\bsm-\bsw|!\bigg)
 \bsb^{\bsm}\, \calS(\bsnu,\bsm) \\
 &= 4\,C^2\,(2\pi)^{|\bsnu|} \sum_{\bsm\le\bsnu}
 (|\bsm|+1)!\, \bsb^{\bsm}\, \calS(\bsnu,\bsm),
\end{align*}
where in the penultimate step we used the identity
\cite[Lemma~A.3]{HHKKS24}: for any $\bsnu\in\indx$ and any arbitrary
sequences $(\bbA_\bsnu)_{\bsnu\in\indx}$ and
$(\bbB_\bsnu)_{\bsnu\in\indx}$ of real numbers, we have
\begin{align*}
 &\sum_{\bsm\le\bsnu} \binom{\bsnu}{\bsm}
 \bigg(\sum_{\bsw\le\bsm} \bbA_\bsw \,\calS(\bsm,\bsw)\bigg)\,
 \bigg(\sum_{\bsu\le\bsnu-\bsm} \bbB_\bsnu\,\calS(\bsnu-\bsm,\bsu)\bigg) \\
 &= \sum_{\bsm\le\bsnu} \bigg(\sum_{\bsw\le\bsm} \binom{\bsm}{\bsw}
 \bbA_\bsw\,\bbB_{\bsm-\bsw}\bigg)\, \calS(\bsnu,\bsm);
\end{align*}
and in the final step we used \eqref{eq:1sum} again. This yields the norm
bound \eqref{eq:final-per} and \eqref{eq:final-K} as required. This completes the proof.
\end{proofof}

\section{Conclusion} \label{sec:conc}

We obtain new explicit regularity bounds (Theorem~\ref{thm:der}) for both the standard non-periodic DNN~\eqref{eq:DNN-np} and our newly proposed periodic DNN \eqref{eq:DNN-per}. 
The bounds depend on the neural
network parameters as well as the choice of smooth activation function.
We see that three important factors need to be under control:
the columns of the matrix $W_0$ need to decay (cf.~the sequence $\beta_j$ \eqref{eq:beta}),
the matrices $W_1, \ldots,W_{L-1}$ need to be bounded (cf.~the sequence $R_\ell$ \eqref{eq:R}), and
the derivatives of the activation function need to be bounded (cf.~the sequence $A_n$ \eqref{eq:sigma}).
We also state regularity bounds in the form that covers three commonly used activation functions: sigmoid, tanh, and swish (Theorem~\ref{thm:reg}).

We provide bounds of the mixed-derivative-based norms of the DNNs (Theorem~\ref{thm:norm}).
By imposing restrictions on the network parameters to
match the regularity features of the target functions (Theorem~\ref{thm:diff}), we prove that DNNs
with $N$ tailor-constructed lattice training points can achieve the
generalization error or $L_2$ approximation error bound ${\tt tol} +
\calO(N^{-r/2})$, where ${\tt tol}\in (0,1)$ is the tolerance achieved by
the training error in practice, and $r$ is related in different ways depending on the settings (Theorem~\ref{thm:err}) to the summability exponent $p^*$ of a sequence that characterises the decay of the input variables in the target functions, and with the implied constant
independent of the dimensionality of the input data. 

We restrict the network parameters during training by adding tailored regularization terms.
For an algebraic equation mimicking the parametric PDE problems, our simple numerical experiments with the periodic DNNs and the sigmoid activation function for a periodic algebraic target function show that
the DNNs trained with tailored regularization perform significantly better than the standard $\ell_2$ regularization. More detailed experiments, including with other activation functions, are considered in our survey article \cite{KKNS25b}. 

We stress that the explicit parametric regularity bounds for the DNNs obtained in this paper are completely general and can be applied under other theoretical settings in conjunction with other methods or scenarios. Further experiments on the actual parametric PDE problems, or with different QMC theoretical settings, or with other training strategies and network architectures are open.

\section*{Acknowledgements}

We acknowledge the financial support from the Australian Research Council (DP240100769) and the Research Foundation Flanders (FWO G091920N).
We are grateful to Ben Adcock and Peter Kritzer for valuable suggestions.

\ifdefined\arxivstyle
 \newpage
 \subsection*{Authors' addresses} 

 \textbf{Alexander Keller} \\
 NVIDIA, Berlin, Germany. Email: akeller@nvidia.com 

 \smallskip\noindent\textbf{Frances Y. Kuo} \\ 
 School of Mathematics and Statistics, UNSW Sydney, Sydney, Australia. Email: f.kuo@unsw.edu.au 

 \smallskip\noindent\textbf{Dirk Nuyens} \\ 
 Department of Computer Science, KU Leuven, Leuven, Belgium. Email: dirk.nuyens@kuleuven.be 

 \smallskip\noindent\textbf{Ian H. Sloan} \\ 
 School of Mathematics and Statistics, UNSW Sydney, Sydney, Australia. Email: i.sloan@unsw.edu.au
\fi

\end{document}